# Cultural Endowment as Collective Improvisation: subjectivity and digital infinity

Victor Peterson II

Philosophically, a repertoire of signifying practices as constitutive of a cultural endowment was said to be ambiguous or unworthy of pursuit. Currently, a unique capacity of the mind is considered to be its ability to produce a digital infinity. The infinity produced, here an operation expressing subjectivity, follows a simple principle according to which a limited set of means, here functions, are utilized to produce an infinite range of potentially meaningful expressions. It is from this concept that I propose a theory of subjectivity and the endowment from which it expresses a self in the world(s) it participates. In particular, I make a case for the subjectivity of blackness. I treat subjectivity as an operation in order to address problems with Identity theory, Afro-Pessimism, and to formalize an analysis of blackness despite the onto-epistemological commitments of racialized systems of categorization. In sum, subjectivity will be characterized as poetic computation.

> . . .us colored folk are like branches without roots,
> and that makes for things come round in queer ways.
> –Zora Neale Hurston

In 1963, Amiri Baraka, then LeRoi Jones, stated that modes of expression are modes of being in the world. In this vein, he said the mode of expression of blackness is exhibited by the blues. As blues people, the finite means of that repertoire were put to infinite use and become embedded in almost every musical form in the Americas (Peterson II, 2018; Rapport, 2014; Johnson, 1925) and rebounded to areas globally, setting the conditions for a black form of life to emerge. This essay will take this articulation seriously. In his work on improvisation syntax and scat Brent Hayes Edwards (2002) stated that scat began with a fall; a subject's musical expression articulated outside of strict adherence to notation. Thus, blackness is shown to be able to articulate a self and maintain subject continuity outside of strict categorization made interchangeable with one identity or another. Here lies where we can make a start towards a theory of subjectivity that averts the failings of identity theory. Subjectivity is the projection of this innate creative capacity, thus the only innateness to be had is an operation, not an essential thing as identity purports. An operation that constructs one's own and, therefore, can impose on others an identifier such that the use, the function of that name entails that subjectivity is primary, not the system of categorization it just so happens to become a part.

The concept of a cultural endowment represents the lexicon from which a form of life is constructed. What is articulated from this endowment has been central to Black study since the 19$^{th}$c. The development of this concept is integral to positing a continuity of expression whose operation is Black subjectivity. This endowment has been posited both as an affinity of signifying practices (Gates, 1988; Hurston, 1934) and as a coalition constructed by way of one's mode of being in the world. (Baraka, 1963) Modes of expression are the means of articulation from endowments. As a mode of expression, blackness is considered a capacity, not a thing (Spencer, 2014) its operation indicative of subjectivity. The method with which it constructs the form of life it expresses can be shown actively constructing that world. (Gödel, 1941; Jones/Baraka, 1963; Wittgenstein, 1922) The problem of blackness is a problem for thought (Chandler, 2013), thus set-theoretical construction, as a way of manipulating the objects of thought, (Cantor, 1874) give us an entry point for the study of the functions whose operation expresses the concept of subjectivity.

We start with modes of expression. For Jones/Baraka, modes of expression are the expression of a people. As such, in accordance with Wittgenstein, the expression is of a form of life. Blackness from here on is not merely the predicament of an absolutist determination of a group, but an affinity obtained between modes of expression indicative of the active construction of a form of life. This is in excess, not just in reaction to an imposed predicament of strict categorization—whether racial, national, or geographic. Blackness will be considered as a disjoint union of multiple signifying practices, indicating cultural and ethnic antecedents that are not centrally placed, but based on that affinity between the operations constructing them. Expressing and, therefore, forming the endowment I wish to illustrate. As a faculty of subjectivity, the formulation of this capacity is the closest we get to making an account of subjectivity as functional, indicated by the relation obtained between input and output claims. (Chomsky, 2015) This, then, composes a relation between subject and context because the capacities of subjectivity are the means of context construction.

We start here in opposition to a presupposed categorical description. To commit the latter would make subjectivity a byproduct of identity without accounting for what identity actually identifies. Ultimately, this leaves us with the question, what constructed that identity in the first place? And if by some other, by what means? And if by some context, how was the constructed context determined? With that being said, there are three factors that must be accounted for in order to come up with an adequate explanation of what a cultural endowment is with regard to the activity of the subject constructing it.

The first factor deals with the data available. Theorist Marvin Minsky (1967) once constructed two-state Turing Machines allowed to run utilizing two symbols to code input and output.[1] Most produced nothing. However, from no one thing, and limited external stipulations, a recursive successor-function emerged, creating different orders of data from the finite means available. In terms of modes of expression, the successor-function modeled merge, the primary operation for the language faculty. Language is not what comes out of one's mouth, but the mechanism for creating thought. From this, we arrive at a natural starting point to make an account of subjectivity. Stephen Wolfram (1991) theorized a two state, three symbol machine proven universal by Alex Smith. (2007)[2] The two state, three symbol machine—A,B states with symbols 0,1,2—revealed that decidable problems are only such after the construction of a reference frame. A yes/no determination is relative to the context of application, the model of the data constructed by some such operation. Considered ordinals, "1" determines the name of context, $\{0\}=1$, from an initial state, 0; and 2, as the set $\{0,\{0\}\}=\{0,1\}$, constructs a relation between input and output. The function constructed is applicable independently of any stipulated norm. As such, the emergent operation is recursive. The finding from the difference that emerged between [2,2] and [2,3] models showed that causation is not identical with correlation between data and the categorical system brought to the data set. If there are only two symbols the context is assumed, therefore indeterminate. By adding a third symbol, a frame of reference indicates the construction of the context in which this operation can be discerned by its activity. So we can understand a normative concept explained by non-normative means, applicable without the norm determined beforehand. The ability to take two objects and make a third without losing the former, indicates the initial state of subjectivity is null when considered from a position external to context; and yet, as it is a function, its zero may be null, but its content is this internal creative capacity constructing the context it inhabits. Constructability is easily shown by Hasse diagrams with the subject occupying the null position, lines as functions, and nodes as constructs from those functions. From that diagram, the shape of the subject's reality as the relation between functions, its mode of expression as its way of being in the world it constructs, can be illustrated. From this we state that a marker of subjectivity is the ability to make infinite use of finite means.

Now our second factor. Subjectivity is expressed by a recursive operation that creates objects appropriate to data but not necessarily caused by the data, therefore not random. With this understanding in tow, our third factor brings us to an explication of what a cultural endowment entails. If it can be shown that this individual capacity holds, then what would be the case of a collective capacity? As a successor-function can be elucidated by set-theoretical construction, a composite function based on the product of this individual capacity is likely. In this way, the subject making use of finite means can produce identities appropriate to context but not necessarily caused by that context. This aids us in understanding why it may be the case that an identity has not been experienced before. Taking into account that subjectivity is an operation, these identities are not seen as random or arbitrary but are the result of a particular operation whose function indicates the modality in which the entities of this particular endowment are expressed. This is the case because those subjects actively construct that context.

First we will show how subjectivity operates. By that operation it constructs the contexts, the realities, in which it inhabits. Then we will show how the accumulated contexts, realities, for the worlds a subject inhabits, are constructed by that same or other subjects, prove an endowment from which that operation pulls in order to articulate identities. Identity is how that subject is in that world, appropriate to each context but not caused by any one context. Therefore, the subject's articulations are not random because the endowment is transitive across those worlds. In this way, particular input place limits on the type of identity produced from the set of all possible. This subjective-operation produces different output

---

[1] *Computation: Finite and Infinite Machines—*see their [2,2] machine.

[2] *New Kind of Science*, p. 709—see their [2,3] machine.

dependent upon the context in which it is situated while actively constructing it. However, as all inputs cause certain things to emerge, all identity constructs indicate the conditions of their emergence. As conditions, the context indicated by the material available to subjects and the relation they obtain with each other due to the function expressed by this operation, demonstrates that changes in the context also cause changes in the domain of their endowments. As subjectivity is indicated by the activity of identity construction and as identity implies a context in which its assertion is applicable, this indicates that the subject-operation and the domain from which it pulls to articulate a self is present before input. The type of entity it is in and only in that context, how it is identified, is what is context dependent. (see Appiah, 1985)

II.

Now for the recursive operation expressing subjectivity. The construction of contexts indicates the activity of such an entity. We treat this capacity as an ordinal. The null, which assumes no one thing, can be shown via some operation to express any thing. Therefore, null can be found as the basis for every context. We take the von Neumann ordinals based on set-theoretic construction as a model of this operation. Starting with an object of no scope, null assumptions—e.g. a name, axiom, object of knowledge—we proceed in the following way,

0 = "null"
{0}= "1"
{0, {0}}= "2"
{0, {0}, {0, {0}}} = "3"

By this method we can conceive of a relation ≤, each ordinal constructed prior, is housed in its successor; each successor constructed from those prior. We can assume 0{0}=1 without changing its being 1. Thus, "1" is the creation of a context by the very function of this subject-operation. "2" relates subject in or to context, for {0, {0}}={0, 1}. Identity is substantiated within frame. The disjoint union of what is the case with what that assertion is not comprises the domain of all possible articulations. Identity, then, is not to some thing external to the operation but refers to the context of application, the context in which "1" may or may not be asserted. Constructs greater than 2 relate contexts with others, relation between contexts construct states of affairs. Identities are discernible based on their obtaining a functional role in and only in the context of their construction. This operation, i.e. subject, is transfinite. As the content of null, the zero of all possible functions, the subject, is both a part of each context by some identity and yet apart from each context by some other. Linking the two is the operation constructing both. The same subject can obtain multiple identities tied to the particular contexts in which they are applicable. The subject is a member of that context if it can be shown as participating in its construction. Otherwise, the identity asserted is inapplicable. As such, the capacities of subjectivity displayed by its activity is able to change contexts because it actively constructs that context with each assertion of an identity therein. The subject is outside of the context of identity. Identities relative to each other compose frames of categorization which can be shown to be constructed. In and of itself, the subject's determinant is null. The function of identity is vacuous without a theory of the operation embedded in its assertion. As a null-state, subjectivity is no one thing. As no one thing and therefore the state from which any thing is constructed, a unity between contexts can be proposed.

Null is at once a part of and apart from every context. Adding or taking away null from context would not change that context's identity. However the context must be constructed from some state. We determine a function for subjectivity from s(0)=0, our constant, i.e. the function itself. From here, we state that there is a primitive recursive function, the "function" of a function, S, such that S(n, s(n)). This, I feel, assumes the least amount of stipulations. What is known by name has no scope as names are axiomatic in the contexts they construct. It only has scope if and when applied. A name only has significance, then, in and only in the context of application.

From this simple operation more complex functions can be discerned, not by stipulation, but by structure. Addition—the union of the elements of two or more sets—and multiplication—the summation of the pairs of elements from two or more sets—are possible after the formulation above. These secondary and tertiary operations from which identity, aggregate of entities under a label for the former, comparison between labels for the latter via a theory of categorization, i.e. quantification of groupings or pairs, makes sense. Identity theory straight off seems implausible because no one to one correlations between individual and category obtain by feature alone. Subjects can be part of multiple categories when features considered separately gain that same subject access to different category.

If $\{0\}=1$ designates a collection of objects/features out of all possible, each ordinal is a determination of that designation. A set of groupings does not always entail the exhaustion of features from each collection entering that operation. Automatically we can rule out Afro pessimism, a theory that assumes Black identity is a cultural death. Either universal identification of blackness as interchangeable with absolute nothing is arbitrary—for names have no scope (Barcan-Marcus, 1961; Kripke, 1972) only the extension of their use or application—or that theory captures no one thing. Both result in blackness' cultural endowment being non-empty. Blackness is not a thing but an operation. As such, this operation populates every context in which blackness participates. We will discuss Afro-Pessimism in further detail below.

It is from this that we imagine that the objects considered knowledge within a particular context are assumed. As merely labels, tags affixed in order to define a function, they have no scope in and of themselves. They are axiomatic within reference frame because they cannot be proven within that frame. They become the basis upon which that frame is constructed. As such, a certain value is attributed to a name. Names are not a necessary constituent of what they label. Otherwise, we would get into a Russellian, if not a Gödelian, paradox: a set whose members are not self-members; or akin to the undecidability of a self-reflexive statement about "a", the system in which that statement is formed. Identities, then, only have a function once applied in such a way that determines a reference frame or context-model of one's state of affairs. One based on the material employed to construct that model. If contradictions arise based on the method those names were employed in order to construct that context, i.e. one's reality because of the consequences the assertion of its identity therein entails, then the model is seen as inconsistent. However, this inconsistency is stored as an axiom in the overall system, becoming that which cannot be proven in that system. Thus, null's unification of all possible contexts is valid based on a sort of Boolean ring of idempotent objects. Here 0,1 decidable in 2. "2" the context or frame as the set $\{0,1\}$, in which 2 is expressed by the functional content of that context, but is not represented in it because it is the context itself, determining the discernible relation between 0,1.

This unification can be shown in the following way. We can treat a label of no one thing in particular as a constant. Thus a of 0, $a(0)$, represents a name as yet to be applied. We know that $a(0)=0$ defines the zero of that name-function as a constant. This formulation can describe the operation expressing this concept. Multiplication can be derived from addition using primitive recursion. Disjoint union, "or," is represented by '+'; and the pair-wise aggregate of the features of each set to which a is applied by exponentiation. Since $a+a$, "0 or 0," in accordance with Boolean methods, equals $(a+a)^2$ then, $a^2+a^2+a^2\ldots=a+a+a\ldots$ Thus, $a+a=0$, for a is just a. With no stipulation, the operation of a as a being constant indicates an operation over null expressing no one thing as any thing. Therefore, a and a, a or a, etc. just is a. We have said no one thing about a which, therefore, is not related or framed by any other thing. For any x, we obtain the same result. If $x+x=0$, then this disjoint union tells us that x and not-x is null as well as for the conjunction, multiplication, of both. So if the ring is characterized by a power of two, we can prove the commutativity of this operation over null, thus relating any one context constructed from that domain. Commutative of subject's over contexts, not the identities asserted within those context, i.e. rings. For x,y with y=not-x. As such, xy=yx for xy+(xy+yx)=xy asserts something about x and y, namely that they are in a relation constructed by some operation.

For any string that comes to an odd power, we discern two or more entities make up the string. We suppose that we have x+y just in case we discern that one of the x's of a string $x^3$ is not-x because (x+not-x) could be (x+(-x)). If no one thing is said of x, all x is unified, i.e. an element, in null. Thus, a ring, context, is proven for x and a different ring for y. From that additional "odd" entity added to an already proven ring, another context can be constructed, including the one previously formed, from which a ring can be proven for that additional entity. We can always form another ring. If contexts are ordinally constructed, chains of context citation are demonstrable, thus plausible, for any assertion because an assertion implies the context in which it was made. If not, then the assertion is incongruent with this chain. The assertion is either of another or starts a new chain from null; or this is its only context, ending where it began because empty or self-refuting. Shown as a commutative property, wherever one is in a chain of contexts, an assertion can refer to a valid predecessor off of which it was constructed and, therefore, is applicable in this context. One can pick up the function of that assertion if one comes to be in a relation with it as constructive of a shared state of affairs. This lays the ground work for how endowments both evolve and are passed down.

What is determined not-x gives us a non-determinate y such that y is the frame in which the function of x is expressed. Given a parameter y such that the function of x is such that, $f(x\sim x'| y)$, we do not have to quantify in that frame. We are able to ascertain the function constructing the scope in which an operation of x is understood relative to y. What is determinate regarding on operation of x in respect to no one thing, $f(0)$, means that no change in x can be determined on its own. That function, however, can be understood as an object not present in that ring but represented in the change between antecedent and emergent conditions. The content of the assertion "x" is the function of that expression, not x itself. For $f(x,x)$ expresses no one thing about x save for that function without application. However, as the zero of this function, this represents the function itself of which a change in x indicates its applicability within parameters set by what is not-x, thus y. We will show the import of this in our penultimate section. We know that this function is likely because in order to ascertain a change in x given x, we must hold what is not-x constant, thus y given the change observed, $f(y|x\sim x')$. There must be a limit to this frame otherwise change would be indiscernible. The likelihood of a function is determined by considering changes in observed data as a function of the parameters, the limits of what constitutes what is observed, thus indicating something producing those contents, of that domain. As such we see the function as the content of assertion made within the context set by y.

The construction of a reference frame is based on the context of the appropriate application of the terms of knowledge. According to P.F. Strawson (1959), the objects of one's knowledge must be discernible and re-identifiable within the contexts applicable. Applicability is such because these are the terms with which that context was constructed. A context of that context, that terms use in a successor-context, is valid if and only if the one prior is consistent—see Löwenheim-Skolem's theorems and Turing's Ordinal Logics. The subject's operation over these terms of identity construction, their being put to use, indicates the existence of a repertoire of functions, of signifying practices that comprise a cultural endowment. Each function's constant, its zero, e.g. $f(0)=F$, represents a term of that repertoire's body of knowledge prior to application. The move from no one thing (null) to some thing (model) is based on the successor function above. Only in this way do we fulfill criteria for an account of the formation of an endowment. In this way, we obtain a plausible foundation to explain what we know (axioms) and the method of constructing models from that set (articulation).

In what follows, we will consider cultural endowments and provide proof of its existence based on the conception of a workspace. We will develop this notion within the concept of a lattice representing the topology of the endowment whose structure is conceived as it is constructed. We will work through the concept of articulation theory—the study of the conditions in which relations of subordination and dominance emerge—in line with Stuart Hall's (1980) request for a formal account of this process. An account for the conditions in which contexts are constructed and interact is made to provide the space in which identities emerge due to the operation of subjectivity. Hall's request was made in lieu of merely documenting or making a statistical analysis of the effects of this process, leaving the conditions which reproduce these relations unquestioned. It is well known in statistical study that correlation with the categories organizing the after-effects of this process, the events experienced, does not entail a theory of causation, only a description of what we already know. Subjectivity is expressed by the successor-function modeling a merge-operation that takes from a lexicon objects of thought in order to articulate an identity appropriate to context but not necessarily caused by context. Therefore, subject-articulation is not random. The same subject can maintain its consistency while projecting multiple identities.

An articulatory-operation takes two objects defined over a workspace producing another. The object articulated is the set of those objects that were once separate or one part of the other prior to that operation. This does not change the objects themselves. Based on the successor-function, this operation represents the creation of the successive construction of well-ordered sets. We obtain with each successor a context (w), in this case a workspace, an object (n), and that objects relation to/in that context (r):

$<w,n,r(n)>$

Formally, we arrive at the notation of a vector. As such, we can map the articulation of subjectivity as having direction and magnitude. Vector notation allows us to conceive of subjectivity as a process without having to make that process interchangeable with some other thing as the object of its identity. If this was done, we would lose the relation and concept we wished to elucidate in the first place. Our analysis would defeat itself on its own terms. In accordance with vector notation, the determinate, if subjectivity is our base, i.e. constant, is: $w-n+r=0$. Thus, 0 is any n, therefore no one n in particular. The basis for constructing models of our state of affairs is represented by the formula: $n=w+r$. This latter formulation represents a context and the means of that context's construction. This, then, is a formula for discerning between subject and identity and reveals that identity has no scope until asserted within the context in which it is used.

The method of constructing workspaces can be shown in the following way. When n=0, its assertion implies its context of application. Thus,

$n=0$
$\{0\}=w$
$\{0,\{0\}\}=r(n)$
$\{0,\{0\},\{0,\{0\}\}\}= <w,n,r(n)>$

where r is recursive. When n=w, we can assume $n=w+r(n)$. As the operation that models the articulation process, this is the union of a context and the means of that context's construction.

If we take two contexts $w_1=a$ and $w_2=b$, we can construct an *R*-function such that:

$R(a,b,w') \rightarrow [\{a,b\}=c]$

where w' is the workspace implied by the operation R on a,b. In this way, [{a,b},a,b] in total is equal to $w_0$ as null, the set of no one thing. As the items available to the operation, {a,b} represents the workspace in which "c" emerges. The items a,b are present but not represented in the successor context, [{a,b}]. However, because we can show how c was constructed we understand that both "c" and a,b as features of c exist, but in different senses. As a,b is not represented in w', if the subject inhabits w', its operation on w' does not have access to a,b; otherwise, contradictions arise. To manipulate the internal structure of an entity from an external stance writes that entity out of existence. The items a,b and their relation to each other expresses the subjectivity of what is subsequently identified as c by that internal relationship. Thus, an operation which takes up "c" as an identified object within the workspace stops at the subjectivity of c, the operation relating a,b, as it only takes the identity "c" as the object of subsequent operations. However, the subject is both a part of and remains apart from the contexts in which it inhabits. Each articulation creates a line of citation $w_1, w_1', w_1''...$ and $w_2, w_2', w_2''...$ which are disjointly and exhaustively valid. As such, the operation constructing them makes a non-empty null set, $\{0\}=W_0$. Their summation by subscript alone are conjointly true. So $w_1 \neq w_2$ but if $w_1$ is valid, so is $w_1'$.

This does not violate workspace conditions following Chomsky. (1995:243, 2017, et. al)[3] Workspace construction is nested so that the scope of the articulatory-operation is well-defined and determined by what state it has reached as per the successor-function.[i] So w=[a,b] such that w=0 or null, and w'=[{a,b}]=[{a,b},[w]]=[{a,b},[0]] such that w'=1 for null=0=everything, 1={0}, 2={0,{0}}={0,1}, etc. Thus, workspace construction follows the construction of ordinals.

By illustrating the operation indicating the subject's articulation of identity brings with it the context in which that identity is appropriate as the expression of that context's subject. We will call this a reference frame. Frames index the context of the assertion of identity as a function of nominal application. As an individual exists in finite time, this limit is placed on the state that subject acquires from the infinitely many possible, making the capacity recursive. The states attained by a subject over the course of its form of life are unique to that individual due to the context in which it participates and/or is introduced. Workspaces represent the reference frames in which the objects available to articulation are restricted. Restriction is required, otherwise the operation would produce contradictions or its scope would remain null because whatever is produced would have no one context in which it may be asserted. The output constructed, here the identity, would apply to every context without any hope of determining any one individual. Producing no one thing in particular, or one thing applicable everywhere thus any where, makes reference impossible. The use of a name is only applicable in so far as we know where that name does not apply, providing a non-normative explanation with the least assumptions or stipulations in order to explain what is retroactively considered a normative concept of reference. Structure provides what is needed for semantic interpretations (Chomsky, 1995) for if we were to rely on stipulation, we would need to provide justification for those stipulations which would not have a method of verifying or justifying themselves. This process avoids those shortcomings.

As the scope of functions are limited, a requirement for producing useful articulations indexing the contexts in which those subject's participate while maintaining subject-continuity, the r-function can be shown to evolve based on past operations. The subject's "stock" of functions ready for application expands with each "function" it determines by way of the context it produces. We can account for both the formation of w as well as how w and w-successors are the conditions for alternative w's. As w is

---

[3] Follows internal merge, here articulation from an endowment, rules:

$f(a,b) \rightarrow [\{a,b\}=c]$ which implies that, $[f_a, f_b, \{a,b\}=c]=[0,0,\{a,b\}=c]=[0,c]$

where {0} represents the zeroed functions of c's endowment and the members of the workspace entailed by the method of constructing that space. The workspace for the subject whose operation expresses c, if $f_a$= "a" as the zero of the function constructing the context in which a was a member, then looks more like: ["a","b",c] where c={a,b}, $f(c)=[c]$, and $fc(0)=$"c" for that label has yet to be applied to/in any one way/object; thus constructing the set of objects considered/labeled "c". Chomsky takes care to emphasize that no thing is lost after this initial merge. Of course, 0 added to some object does not change that object; just as well that object is such because constructed from a function whose initial state is zero. Internal merge, by definition, occurs when b<a, where c takes up $f_a$; or a<b, where c takes up $f_b$. If the operation under c takes up $f_b$, this represents an extension of b as c; if that subject takes up $fa$, no extension of c but a deepening of its endowment as the stock of prior participation- -active or emerging from--regarding context construction thus being introduced to and possibly acquiring those functions: $f_a$={"a"} and $f_b$={a,{a,b}} the function of whose union constitutes what is considered "c"={$f_c$}, the set of zeroed functions previously acquired by c, an effective representation of the subject c's endowment prior to articulating any one thing that is "c"—see discussion on prefunctions.

constructed by r(n), the use of n implies w because it is the material from which one constructs an initial w—e.g. the move from n=0 to {n}. The result, *r* indexes n to w by use in relief to all options available within the union of prior workspaces.

A demonstration of this can be shown by intuitionist logics. (Gödel, 1958; Heyting, 1930; proof by Kripke, 1965) By induction, if we suppose a, then by its assertion we get the context {a}. From this, the construct of this context is understood in relief of others which are not yet barred from doing the same. So within that set of others in which {a} is now a member, {{a}}={a,{a,b}} which we know from Gödel orders a partition of null, <a,b>, where a is less than b. This can be interpreted as a being function of b for <a,b>=<a<b>>. For example, if 3,2,1 is well ordered, then any partition of it is in that order. So, if <3,2,1>, then <3<2,1>>. The domain of b is a. For, a is both part of and apart from b as determined by an ordered-relation in respect to, so pairing a with b. We showed the construction of a workspace and an accessibility relation from structure alone. As well as the emergence of the r-function described above. With r(n) being recursive, we can illustrate the construct of any context as well as how they interact with each other.

Suppose from prior articulations we produce c. We have a function built using the conditional, if/then symbolized by '$\rightarrow$' and '$\times$' as conjunction, such that,

*f*: a$\rightarrow$b

If this operation is transitive, we map c onto this operation,

c(*f*): c(a)$\rightarrow$c(b)

which is to say, c$\rightarrow$(a$\rightarrow$b). From this we can say if there is a function mapping b to a, then if *f*(a,b) then *f*(a)$\rightarrow$*f*(b) for our notion of if a,b implies that a is included in b, the move from no one thing to any one thing as a morphism of that domain. In this project, we call this an articulation from that domain: realities as various context qua models arranged from the objects of the actual world; the set of all objects in the actual world being the set of no one thing, thus null for any one reality.[4] Then if c$\rightarrow$a and c$\rightarrow$b, then if a and b, infer c. We obtain modus ponens for if c given a and b, then if c then if a and b then c. We also see how c emerges from a,b which is both a part of c and apart from c as its conditions. Thus, the endowment that we speak of is a repertoire of *f* and *f*-successors which are not present but represented by the emergence of c within the state of affairs composed. The context of "c" is the well-defined function, the operation over a,b, determining the domain, the context of its use, as b to a. The function *f* becomes a map of the conditions a,b and what is articulated from those conditions is continuous just so long as *f*(b) is included in c. We will see this in regards to our discussion of internal merge and rearticulation below.

---

[4] Let A be the actual, i.e. no one object in the actual world, and let F be a particular arrangement of the objects whose relation express the function of a model qua context in that world. F is included or equal to H. F follows these properties:

1. F is valid
2. If x,y are included in F, then the conjunction of x,y included in F,
3. If x is in F, y in A, and x≤y, then y is in F, expressing that null is both a part of remains apart from whatever is articulated as what is articulated is constructed and expressed by an operation over null.

Therefore, given any subset S of A, there is a smallest contextual reality containing S. We say the context is articulated—for Heyting, generated—by operations in A by S. If S is in null, F={1}. Otherwise, F is equal to the set x in A such that there exist y in S with y≤x. If A follows our Heyting quotient and equivalence rules, and F is a context on A, we determine a function ~ on A as follows: we write x ~ y whenever x$\rightarrow$y and y$\rightarrow$x both belong to F. Then ~ is an equivalence relation such that an operation in null is equivalent to the function of an articulated identity modeled in the context constructed. Thus, we state that a model is such that F/~, whose elements are represent [f]. Functions utilized in the reality labeled by the functional content of F, [F]—as [$f_0,f_1,f_2,…$]=[$f_0,[f_1,f_2,…]$]=[F,[$f_1,f_2,…$]] and F=$f_0$ is the label whose application is a function of all objects in F—are determined by an equivalence relation on operations in null conceived by set-theoretical construction. So, from f$\rightarrow$e in null we get if [F$\rightarrow$E] such that A(F,E), where F,E are related by the function determined by operations in A. We will go into this in more detail below.

According to intuitionist logics, a→b is $b^a$, for a→b=a≤b. We have shown above how the concept ≤ can be defined by set-theoretical construction. So, if we have values 0, ½, and 1, then we have a→b so long as a≠0, otherwise b is any one thing. If $b^a$, then the set of everything is no one thing for b, just so long as a. If a=0, then a must be the set of everything for no one thing as yet has been determined. All subsequent articulations follow from null, the set of no one thing framed by b. Thus c→(a→b) also shows c(a×b), inferring an operation over a,b producing c, all couples, possible arrangements, of features of a,b articulate the function that is the content of c. Boolean conjunction expresses multiplication. So, if not-(a×b) then not-c because there is no a nor b. Therefore, $c^{(a \times b)}$, c a subset of $b^a$, whose extension is up to where c just is $b^a$. So, there is a domain of c, even if null, that is $b^a$ given the articulation of c. In other words, there is a c within the frame constituted by (a×b) the function of which, $f(a,b)$, is the zero of c, $f(f(a,b))=f(c,0)$. This is the functional content of the context of c's assertion, $\{f(c)\}$, before any one application of "c" producing the entities considered c.

How does this apply to articulation theory? Say that a=<w,n,r>. From our method of construction, n, then {n}, then {n,{r,n}}. We do the same for b. The context of a's assertion, in accordance with Gödel, would be a=<w<n,r>>, where w is expressed by a function of some name whose operation indicates the means of w's construction. The workspace w, then is determined as the group (a,b), so (w(a,b)). An operation on that group would entail,

R(a,b,w)= (a×b)→w'=[c={a,b},a,b].

Let's say that c enters into an operation with a. The constituents of a do not enter into the operation, only what is represented on the workspace {a,b}, the union of "a","b". Their names alone make the zero of the function "c", positing null as both the conditions from which c emerges and part of c. So the operation is ({a}×c), c and the name of a, which is not c×[a| <w,n,r>]. Otherwise, c would change its internal structure making c≠c, canceling itself out.

As {a}=n, for "a"={$f(a)$}, the object that enters into the articulatory process, the identity of a, is separate from the operation constructing a; its subjectivity exhibited by *f* cannot be the features of [a], but are present in the conditions allowing us to discern that workspace.[5] As name, $f(a)$, being the content of the assertion of that name means that *f* produces all possible subsets constructed in that initial context implied by that names assertion: {a}, {a,{a}}, etc. For a's features, the content of its assertion is the function relating those aspects of a, are not represented in the workspace to which c has access. [a] is indicated by the operation constructing the set representing it in total from elements outside the scope of c. In this way, subject continuity can be shown so that a is non-empty. Again, blackness is not vacuous as the Afro-pessimist claims. Subjects can only be labeled, identified, as such. Negation presupposes the existence of the entity it negates. An identifier, "empty," is not a necessary constituent of a for a is only such by way of its mode of expression, the operation indicating its subjectivity, expressing its form of life. The identity is merely a label applied over or within context. Subjectivity cannot be empty for that would negate everything that animates that label. The subject must not be empty in order to make its categorization as "nothing" possible at all.

---

[5] As null is the only set that contains itself, the content of an object in itself, as no other thing, is the function of its expression. The zero of a function is no one thing, thus a member of the set which only contains itself as a member. The set of all functions, thus no one function, uses a function determined from null to label any thing to construct the context it inhabits. The subject itself, the "floor" of that context from which its operation is expressed, is null for it is no one thing to anyone in context. The identity itself, the "ceiling," is the determination of the function of that operation, the least object greater than all that is in that context it constructs. E.g. <$x_0,x_1,x_2$>=<$x_0,<x_1,x_2,...>$>where $x_0=f_x=X$ for embedded in the assertion of X is its function of application. Without arguments, X is a label whose application is appropriately determined such that <X<$x_1,x_2,...$>>As such, it is a label applicable in that context. Our notation, "a"={$f_a$}, represents the operation expressing a's subjectivity. The function of "a" emerges from an operation within null, yielding valid identifications of that subject dependent upon context of assertion. Embedded in that function is the limit of the extension of various articulations of "a" which need not be predetermined. In lieu of using ⌈…⌉, for the function of identity we will use [...], as an equivalence class modeling underlying operations. Brackets note the interval over the arc of the form of life constructed by the operation of the subject whose function is that subjects identity in that interval, i.e. context. For the subject-operation from null, the symbols ⌊…⌋ are replaced by {...} to suggest set-theoretical construction. In this way, we define equivalence classes whereby operations over elements of null represent a function over classes by virtue of Heyting algebra rules. So, ⌊a⌋=⌈a⌉=a. By definition, if x varies over contexts of assertion such that 0≤x≤1, for n<1, if x<n *floor*(n)=x, and for n<x *ceiling*(n)=x. So, {$f(a)$}, the set of all to which "a" applies, displays subjectivity (x) regardless of identity (n), for the subject-operation is the function $f(a,0)$ and the range of identities of a is such that $f(a,x)$ where a≤x and x is the context in which "a" is asserted, indicating the same subject expressed by *f*. It follows, that null is both a part of and apart from every assertion whereby contexts are constructed, and therefore inhabited, by some subject-operation regardless of identity stipulation. Subject capacity outstrips identity stipulations for beyond the ceiling, there is some operation whose expressions cannot be identified. We know this because the fractional-part of this function states that, unless subject and identity are equal, null, then there is always portions of the contextual-floor, that is not accounted for by the ceiling-stipulation within the set whose internal operations express subjectivity. We affirm that there is an *f* at zero for all entities. $f(x_0)$ is the subject within context before identity and as such is included in the set of x. X is the label applied to set of all that is x; all that is in excess is not-x, thus X represents a ceiling to the category of what is considered x. X= $f(x_0)$, representing the subject that has acquired that function who has yet to apply it to, utilize it to express, any one thing, thus the subject is no one thing in and of itself. As such, that function, in itself, is equivalent to the label applied to any one thing, utilized to denote no one x but every x. Reference cannot be uniquely determined save by reference frame, otherwise leading to contradictions embedded in the act of identification.

The function constructing "{a}" is an operation of a, making the content of { } the function of a. The identity "a" does not equal a. If "{a}" is such that [a] is a member, then [a] is a member of {a} by some property. That property is not-[a], thus [a] is not a member of itself but a member of the set denoted by that property. That set is only such by way of its method of construction from all [a], otherwise a would not be a member of itself for it would be identical to that property, making the absurd statement that a is "a" if and only if "a" is not a.

Setting up the workspace example in this way, the concept of reference is taken care of as a superficial descriptor is applied after some constructive operation. If reference was merely the exchange of terms or the interchangeability of an operation for the term it produces, reference loses the relationship we wished to precisify in the first place. A chain is created by the function illustrated above, providing the correct semantic interpretations required for a reference frame by virtue of the syntactic structure of the contexts constructed. In this way, *f*(c) is seen as embedded in the workspace, the endowment expressed by an operation over its antecedents. As such, c is an extension as well as a deepening of that subject, only now identified as "c". The embeddedness of c within the conditions from which it emerges—see the construction of ordinals—entails that the means to rearticulate states of affairs are embedded in each context by virtue of that context's having been constructed in the first place. Any stipulation of universality sows inconsistencies in the system held constant over contexts for, inevitably, there will be a context in which c and not-c are the case by virtue of making that name interchangeable with a context in which it is not applicable because it as its function was not the means of its construction.

Suppose that "c" becomes a function whose domain is a,b and produces a successor c'. We translate a '→' as '≤', conjunction by multiplication, '×', and assign values to expressions that range between 0 and 1.

---

However, functions can be uniquely determined in a way that they become available for application within the appropriate situation because previously used to construct the context that situation cites. The object of thought *f* is the function whose application populates expresses the x in the category set by X

We can assume, then, that $f(f(x_0), f(y_0))=f(0)$ is same subject, $\{f_x, f_y\}$, the set of functions acquired and yet to be applied. This endowment identified in different ways. The "knowledge" of self, the zero of various functions, has the capacity to assert alternate identities so that $x_1$ is not equal to $y_1$. This shows that null is both less than, included, as well as greater than, external, to every identity constructed or identification made. E.g. $\{x\}=x-\lfloor x \rfloor$, then $\lceil x \rceil = f(x)=$"x". From this it follows that there is always a remainder of x after "x" is calculated, implying another identity under which x operates inaccessible to *f*(x). The contextual-construct "me"="Victor," is me based on prior subject-representations of "me." It follows that I am in excess of both those prior representations as the "I" constructing/using them as well as the context constructed/identified as me—unless this is my first cognizant moment, implying the existence of an initial state whose content is a prefunction of which its prompted operation is the subject structuring the initial data set, perceived as no one thing externally, but whose function, after application, is how that subject is identified by others. The function of that identity is the limit of my participation in that context. This does not violate the fact that the context's limits is the representation of "me" and the function of that identity "Victor" such that I, myself, both participate in and yet remain apart from that context by virtue of other identities. We are licensed to say $\{f_a\}$ because the identity of ceiling and floor is a relation that supports the assertion of an identity as the context that subject participates; even though it, itself, cannot be represented in that context because its operation constructs that context.

a→a = a ≤ a

We define a function of a,

$f_a$: $f(0)$→a

for,

0→a=0 ≤ a

$f_a$=$f(a,0)$=A which is the determinate, the zero, of the label for what we come to call as a function of that label's application, "a." If a is our constant, $f_a(0)$, then, any thing follows from no one thing. We proceed with a statement of b,

a→b = a ≤ b

making for $b^a$ such that we have a frame of a in which we can discern the articulation of various entities.

In our Heyting scheme, "if a then b" is equivalent to a less than or equal to b, which is to say b given a, (b/a). If a is greater than 0, conditionalizing over null, the set of all possible, obtains some partition of null as our reference frame determined as the domain $b^a$ for subsequent operations. If a=0, and with the assumption that any thing follows from null, (b/a)=1/0, which leaves the conditional non-determinate. 1/0, however, is the representation for infinity, proving the point that a non-empty null, meaning null for a singular context, for no one thing in any particular context, is the ground from which anything may be articulated. An operation over null, then, is infinitely recursive, therefore generative. As such, if a is greater than 0 or just 0, the disjoint union of (a,not-a), makes an operation over that space a modal operator, possibly-a. This operation is the mode that "a" appears in the world, articulated in the reference frame of b. The domain of null available for subsequent articulations outside of the domain of a, in what is not-a, is (1-a), for which we will see the usefulness below. If c results from some operation in $b^a$, namely (a×b), then c' is possible in the partition (1-a), for, if a then b is a less than or equal to b, then, c and a is less than or equal to b, for b is set as the domain in which c is obtainable. It pays to remember that according to Heyting definitions, if 1 less than or equal to 0 then a entails that a, for null remains in excess of the entity articulated within the domain stipulated. For if a then b entails $b^a$, and c and a is less than or equal to b, then c is less than or equal to a. Thus, a is both part of and apart from the extension of c. The remainder indicates the existence of alternative contexts outside the domain currently determined, possible within and from null.

For now, we state that a=b, so if 1 then 1, as the domain of our world and within which possible realities as contexts are articulated. "If a then b" creates the frame in which all possible worlds can be articulated, this framework states the conditions one moves from no one thing to a context determined as no one thing. Thus c,d, and others can be articulated in that frame, indicative of their own frames—see ordinal logics, c-successors within b—and are valid because able to be substantiated as a b-successor. Each operation within $b^a$ conducts a pairwise grouping with what is articulated basis in b: [c,b'], [c or d,b''], etc.; for, for all articulated in $b^a$ is valid because b is valid, 0 to {0}, 0 to 1. As such, there is always a part of null that lies outside the scope of the contexts articulated.

We model c's articulation,

(a×b) →c=(a×b)≤c

What if c enters into an operation with its antecedent a creating c'. We say that c takes up a's function, $f(a)=0$, the capacity to articulate "a" appropriate to context, as part of its endowment. As c was articulated from those prior conditions; the zero of $f(a)$ in c is plausible. It does not take up a in and of itself. So,

$(c \times a) \to b = (c \times a) \leq b = c \leq a \leq b$

which makes sense, for

$c \to c = c \leq c$
$(c \times a) \leq c$

then if $(c \times a)$, then $(c \times a) \leq c$, a deepening of c's endowment.

As such, if c results from $a \leq b \leq c$ then $c \leq c$. Otherwise, if $1 \leq 0 \leq a$ where not-a=a$\leq 0$, makes (a×not-a)=0. If $(c \times a)$ or $(c \times b) \leq c$ then $c \leq a \leq b \leq c$. If c is greater than a,b, all that is in excess of that operation is $f(c)$ determined by the operation $(a \times b)$. From the less than or equal perspective, what is categorized as "c" in a category in the frame b to the a is what is generated within that frame by that operation. It goes to say a universal description of c captures no one thing about c, for as of yet, that function has yet to be taken up by any operation. Just the same, "c" is only such by its mode of expression, its application to the outputs of the operation that articulated it in the frame b to a. "$f(c)$" is the name of the property of having been produced by the operation b to a, and as such, what is "c" and yet greater than b to a is merely the label of the category "c" which is not a part of b to a itself but by the entities produced by that prior operation that are labeled such. The label cannot stand on its own. This way we can infer functions producing the constituents of states by their operation. Thus the null set is not void but contains functions (partitioned endowments) as the content of every construct by a subject from that set of no one thing. So, "c"=0={$f(c)$}. For c cannot be greater than itself but that valid proposition, $c \leq c$, states what is in excess of the entities under that label and nevertheless exhibit c-ness by the operation utilizing that label, therefore inferring a function, a subject, constructing what we consider c in the first place.

The statement, if "c" then c, is computable, not just an affirmation of the consequent.[6] We handle this through diagonalization. Diagonalization is ultimately about self-reference, so we state if $(a \times b)$ then c, then if c, then if $(a \times b)$, then c. Firstly, a×b can be represented for expository purposes as, $f(a,b)$. "c" is the representation of the function of c. C, a formula whose functional content is "c". We know this formula to be "c"=$f(a,b)$=$f(c,0)$=C for, if $(a \times b)$ then c, the zero, determinate, of $f(c)$ represents a label for c before application.

By Gödel's computability theorem for diagonalization, we know that C=("c"×c), the pairing of every symbol "c" with $f(a,b)$, the existence of at least one c for an n-ary member of the pairwise union of all a,b. In this way, n="c", the diagonal of n implying "c" and c. This is an affirmation of the conditions and the possibility for the emergence of c because "c" is constructible, computable. This is not to say that "c" is necessary to what it labels.

Let $f(\text{"c"},n)$ be the case where n picks out one pairwise constituent of $(a \times b)$. Assume C represents the diagonal whose functional content is *diag*, D the formulaic determinate of C. Let F be the formula whose functional content is the operation that implies y such that F=($f(x,y)$, C(y)). Assume D's equivalence to the formula, x such that x=("F",F) representing the diagonal of F. Let "c","n" be the representations of F,D respectively. By definition, *diag*("F")="D", thus C(c,n) is valid. Furthermore, D is equivalent to the existence of x such that x=c and there is y such that $(f(x,y),C(y))$; which is equivalent to there being a y such that $(f(c,y),C(y))$. Because of the function of C and the validity of $C(c,n)$, this formula is equivalent to $(f(c,n),C(n))$, which is equivalent to $C(n)$. As, F is equivalent to $C("D")$, then D iff and only if $C("D")$, the function of c is constructible, thereby determinable within the context it constructs.

---

[6] A short proof follows this remark but for a fuller exposition see Gödel's 1931 proofs.

In this way, one can utilize the function whose operation is the act of applying "c" without its obtaining any one object. That assertion was made in excess of the appropriate domain because it was not done in the context determined as such due to its being constructed by the operation with which c was expressed. Although we do not pursue a descriptivist theory of names, we still obtain the benefit of Russel's theory of definite descriptions. For the name "c" is not c in and of itself; nevertheless "c" is c for us when we apply it to those objects in the contexts for which the concept c is applicable. With the assertion of "c" implying a context in which it was asserted, that context can be empty save for its functional content, the assertion of that name.

This makes sense as the functions of a,b construct c as the pair-wise aggregate of "function(name-function,context)=assertion": $f(a)=0$, $f(a,1)=a$, $f(a,2)=a'$ etc. and $f(b)=0$, $f(b,1)=b$, etc. such that the function in itself expressing all that is a is 0 for context construction is ordinally based on the set-theoretical construction of numbers—sets and set of sets as determinate of functions and their domains. So, $f(0)=a, b, c$, etc. However, by some operation, $f(c)=f("a","b")$. This should not be controversial as $("a","b")=(f_b(0), f_a(0))=0$ such that $f(c)=f(c,f(f_a,f_b))=f(0,f(0))$ in the domain $b^a$. In this way, "c", articulated in $b^a$, entails that $f(c)=f(c,0)$ by functional-composition.[7] The assertion of that name, for which the context implied was constructed by the use of other nominal-functions makes it such that "c"=$f(c,f(0))$ following primitive recursive principles. Utilizing the set-theoretical construction of functions: a second-order articulation from that null base returns the constant of the domain previously constructed as $f(c)$ in accordance with functional composition, here the evolution of the functions within an endowment constructed from prior articulations, showing that it is possible to pass a subset of the arguments of one function to another function; and, the role of functional projection, here the assertion of identity as the construction of, because it implies, the context in which that identity is asserted. So, if $f(a)=0$, then by set-theoretical construction c has access to a at $f(0)="a"$. The "features" of the endowment of c, the domain over which $f(c)$ runs, includes $(f_a, f_b,)$ as zeroes, which is just how "c" appears as the continuation of the contexts previously constructed and known as "a" or "b". So, $f(c,1)=c$ includes either the zeroed function of a or b.

In articulation terms, $(a\times b)\to c=[\text{not-}(a\times b) \text{ or } c]$ because $(a\times b)\leq c$. Effectively, if $c\to a$, $c\to b$, and if we have c, for c given $(a\times b)$, then $(a\times b)\to c$ which represents our rule of inference: $c\leq((a\times b)\leq c)=c\leq(a\times b)\leq c$. We understand "given" for we "know," have acquired, the zero of the function whose application actively constructs the context in which c emerges. Thus, c given a,b is not identical to a,b given c, however, given what is considered c is only by the function constructing the context of its assertion, then there is a function of a,b.

---

[7] Remember, for any set "a", we understand that a=0, and "a" is constructed $a(0,x_1,x_2,x_3,...)=f(x_1,x_2,x_3,...)$ such that we construct the members of "a" off successors from $f(a)=0$ on. This can be assumed for any set of objects considered a within that reference frame. It follows at a,b,c,… all have an $f$ at their zero for they are only such by the means of constructing the set characterized as such. Therefore, if $a(0,x_1,x_2,x_3,...)=f(x_1,x_2,x_3,...)$ and so for b, where $f(a,0)=a$ and $f(b,0)=b$ then if $f("a,b",0)=c$, then $(f(c), x_1,x_2,x_3,...)=f("a,b",f("a,b",x_1,x_2,x_3...))$ so that $f(c,0)=c$. Identity projection wise, $f(f(a,b),0)=f(f(a,0),f(b,0))=f(c)$. We understand that "c" because we know the conditions in which c emerged as the functions of those conditions. Having acquired that knowledge, zero of those functions, when applied articulates c. And so, we see how functions are objects of thought, which can be determined and therefore available to construct our world but are not physical "things" in that world but nevertheless must be in order for the world to be anything for us at all. Once constructed, the citation of these functions become available for future applications at the zero for a successor-context's construction. Concepts like "memory," spirit, one's endowment, their content is functional, not a stockpile; hence why unique objects of the mind, functions, are only evoked in unique contexts because they were, thereby that individual was, previously involved in the construction the succession of contexts leading to the one we currently occupy. As such, subjects, and the realities they construct, are, in themselves, incomplete a la the quantum physicist Niels Bohr; cannot simultaneously measure their location and speed a la Heisenberg. The subject is a vector. From an epistemological view, all that can be known is its projection, however, the output of a recursive function, that projection cannot stand on its own and is only supported by the operation immanent to it, its subjectivity. Their content, as such, is a collection of determinant functions, recursively constructing their mode of being. The subject occupies the null position at theirs, or when introduced to the line of successors of another's, construct of the web of relations that is their world, creating outward from there and bringing with them the objects they utilized beforehand. Use, the operation, of functions as content of assertion can be uniquely determined and proven based on context and their successors of assertion based on ordinals; whereas reference relies on stipulations without means of justification on their own becoming either more ambiguous with successive application or inconsistent with previous use.

With a≤b and (a×b)≤c, then if we have (a×b), we obtain c. The "c" articulated becomes a function of the set {a, b}, the pair-wise aggregate of the features whose relation expressed a,b, As such, we obtain c(a×b)="c". It remains to be shown how this property of the operation over a,b expresses the function of "c" whose domain is $c^{(a \times b)}$, stating that (a×b) is the horizon of possibility for c and a function of a. The semantics of this articulatory process covers that (a×b) is necessary in the possible world c, which is only the case when (a×b) is possible. Thus, (a×b) becomes a member of null, both part of c, for it is the axiom from which c was constructed, and apart from c, as the features of a or b inclusive are members of the set of no one thing—see null as basis for the axiom of foundation (Zermelo-Fraenkel) and von Neumann-Bernays-Gödel's notion of ordinals based off Peano's axioms of arithmetic.[ii]

A logic set up in this way formalizes a bounded lattice as the workspace One actively constructed by a mechanism mapping the operation of a subject. Within this space, conjunction, (a×b), are the points at which lines meet, the lines between those points represent conditionalization within that space, constructing the notion of membership to a category expressed by the meeting of lines of thought. The system of categorization implied by the internal structure of the space indicates the reference frame shown by the exponent of a particular individual under some name. So, if a→b, then the accessibility parameter to the domain of that workspace is $b^a$, displaying the conditions in which c can emerge.

From this we show how the contexts in which entities emerge are related constructing that state of affairs. Take $f(x_0)=X$, where X is a nominal tag for any one x, for $<x_0,x_1,x_2,...>=<x_0<x_1,x_2,...>>$ in such a way that, $<X<x_1,x_2,...>>$. From this we say that each x is X related to itself, X(x,x), in every iteration of the appropriate context of assertion. X is the reference frame qua context in which the assertion of x as a function of its application within that context obtains. Thus, each assertion of an x is paired, comes with, $x_0=f_x$ as the function of appropriate use of that name, just so long as it is within the domain of x-successors characterized under the context of X and its successors.[8] "X" as a label only applies to objects $x_n$, produced by some operation whose function expresses X-ness by the function of applying X, otherwise the class of objects x is empty for X applies to no one thing. Even if considered purely a predicate, X, on its own, is revealed as only having a function as content; relating, as arguments, an entity and the relation that entity obtains within the context constructed by the assertion of that label. If we consider the function $f$ defined over sets of entities in x,y such that,

$f(x,y):(A \times B) \to C$

the features of those sets as variables x in A and y in B. If we hold that a is fixed in A, the articulatory function becomes,

$f(a,y):B \to C$

which makes that function some feature of C such that $f(a,y)$ is an element of $C^B$ now available in the domain accessible to the repertoire over which the process of articulation runs. If we let a vary, this gives us a map of that domain such that,

$f:A \to C^B$

---

[8] For $x_0=X$, we show $x_{n+1}=X'$ such that nominal assertion $f(x,n)=X^n$ such that all $X^n$ are included in X securing citational consistency.

in which the concept of subject continuity is maintained. In fact, this continuity is unique such that if $f:A\rightarrow C^B$ is the parameter determining the access the articulation-function has within the domain A,B then there is a function $f_a(y):B\rightarrow C$ that indicates the operation expressing the subject as a function of articulation:

$$f(a,b)=(a\times b)\rightarrow c$$

This function exemplifies the property of articulation. In this way $f(A\times B,C)$ is equivalent to $(A,C^B)$. We, then, model the formation of states of affairs A,B,C, as a function of the relation of individual contexts x,y,z.

A diagram of this process can be inferred from the following principles. From if a, then b, we derive a, or not-a, or b as nodes or points within the lattice such that b is connected to a and not-a and where a and not-a are not connected. When a and b are connected, we get (a and b) connected to c; where disconnected we obtain not-(a and b), which is to say not-a or not-b, connected to c. In this case, c connected to b is c to the power of b. This entails that c is such up until not-b, for a is less than b, c is less than a and b, and where there is no b, there cannot be c. Each connection, then, represents a conditional. The direction of conditionalization indicates the consequent from antecedent conditions. We can think of a lattice as a sort of quilt, its patterns are the repertoire that a subject cites and utilizes as a function expressive of itself in whatever environment it participates as an extension of prior articulations.

In sum, there is a function such that we move from the set of no one thing, "a," to some thing, "b." That process sets up the frame in which c emerges. As "c" becomes that frame, for c is less than or equal to c, "c" may determine a constituent within that frame for which that name is not the constituent itself but only the representation that c constructs of either a or b. For if c, then if a,b, then c means that the emergence of c indicates a function whose operation over a,b are the condition of c—see Appiah on conditionals and functionalism.

Ambiguities in the concept of a workspace—as noted by Chomsky, (2017; 2019) et. al—can be accounted for in articulation theory. Determinations erroneously based on identity rather than subject capacity produce the means of their dissolution. We have accepted the limitations of a finite capacity being required for the production of infinite uses. This ambiguity we will call the identity paradox. Labels on their own have no scope save for their application, but what licenses that application is the frame of reference in which the features which animate the function of that label are negated in order to trade that body of features for sake of membership within category constructed within that frame. Identity, as far as we can conceptually determine its function, is expressed by input and output where that relation holds necessarily within context. If $(c\times a)=(c\times b)=$"c" for c becomes the identity of $(a\times b)$ when considered a pair-wise grouping of the functions as features of its endowment, inherited from a,b. The identity determined by the function "c"=c(a×b) expresses the relation between a, b that denotes "c" if and only if $c\neq(a\times b)$, for if we lose either a or b we would no longer have c at all. "c," then, is the function c in the context $b^a$ and is only expressed by those features whose relation within that domain projects that property. The assertion of "c" brings a,b with it. However, (a×b) are not represented in the assertion for (a×b)→c=[not-(a×b) or c]. Identity is only decidable relative to the workspace, the partition of null entering into the computation. Individually, a label can point to any one feature of an individual, however, through the set-theoretical construction of this operation, the function of that label is structurally, contextually, determined.

III.

To recap, endowments are constituted by functions. The functions are the well-determined objects of mechanism for creating thought as they are the only discernible content of an assertion, the application of these functions expressive of subjectivity without reliance on categorical identity for application constructs the contexts in which systems of categorization possibly obtain. We have framed cultural endowments in this way to show that subjectivity is conceived by the relation between contexts in which an identity assertion obtains, the function of its construction considered the fabric upon which the operations of subjectivity can be considered vectors having magnitude and direction. We further develop this notion via constructive analysis under the concept of presets—see Bishop and Bridges, 1985—the possible contexts in which identities are asserted. The function of identity application constructs the contexts in which those identities are appropriate. The concept of endowments under the notion of its objects being prefunctions. The latter perceived functional content of the assertion made exhibiting the "S"-ness of the subjectivity working therein must be in excess of the parameters of that context by virtue of the use of an identity appropriate to that context. (*n*Lab, 2019) As identity assertion implies/refers to the context of its assertion, the possibility of a context to be constructed (preset) implies there is a function present but not yet represented (prefunction) because no context, as of yet, has been constructed. The operation (subjectivity) utilizing these functions (endowment) as the content of its assertion (contexts), is provable before assertion. In this way, functions become the objects, the so called "features," of subjectivity as the mechanism to create thought, hence the subject, is readily ascertained by modes of expression. Operations over these functions create unique, well-ordered, sequences of application, e.g. methods of computation or algorithms, as unique modes of expression. Thus, as the means of expression, these functions represent a repertoire of signifying practices we have called the endowment of the subjectivity expressed, its lexicon. Although of the same endowment, the same subject can express identities appropriate to context but not necessarily caused by context, even if using only what is available within context. Therefore, subjectivity is a creative capacity. We chose this route because functions can be well defined while identity, name to subject, remains ambiguous. A single subject can have multiple identities; a single identity can be utilized by multiple subjects. Subjectivity as a vector, is an operation expressed by its movement in the world(s) it inhabits because constructs those worlds.

Formulating subjectivity as such, subjects can be analyzed in accordance with their movement. We do not have to consider an interval over the arc that is their form of life as definitive of subjectivity in total. It is this movement that expresses the subject relative to the stock of identity formations, understood by their function, previously constructed and made available within particular parameters. Subjects are individual but not mutually exclusive from each other. The operation constructing contexts implies participation in that context. Although intuitionist logics are undecidable for first-order descriptions, we can model the formation of conditions in which subjects emerge by virtue of identity constructs in a way sufficient to our project. This does not necessitate the stipulation of one or another thing obtaining regardless of context, only the formation of those conditions in which these objects emerge and in which identity becomes decidable is required for a theory of the articulation of subjectivity. Subjects have been posited as vectors in order to capture their movement without having to trade a determined interval, whose coordinates are interchangeable with that subject, at the expense of the internal operation articulating it.

With cultural endowments as the domain over which the operations of subjectivity runs, the distribution of an intensional constructive capacity. Regarding racial stipulations, we are talking about the extension of a system of categorization concomitant with the ontological commitments the utilization of that system entails is troubled. Without the former, the latter is empty. I will continue with a discussion of the former for the remainder of this project.

So, "if a then b" is the articulation from no one thing to a complete description of everything—a description not being the entity itself but a property which that entity, by any one of its features, indicates some category as function of that property; gaining them access to another category indicated by having another property—that description can be of anything. The validity of the function of conditions states that it is not the case that a true antecedent and a false consequent obtains, all that is, is b while the zero of the function a becomes the name of that domain. If an identity is universally applicable, it is useless. We conclude that identity does not refer to some extra-mental trait but to its context of application. Each articulation procedure creates a context within those condition. If a, then there is some context A in which a function over the entity a is applicable. We move from no none context, thus every context unified at 0, and which at present is $A_0=0$, to $A=A_1$, $A_1'=A_2$. . .$A_n$. $A_0$ means that no one and therefore any one feature is applicable. No context, as of yet, has been articulated. The successor A' is constructed by taking the universal description of its predecessor, a context-reflexive name unprovable in that context for defined by it, as an axiom for the construction of that successor-context in which that predecessor can be substantiated. Thus, null unifies these properties which are disjointly exhaustive, $(A_0, B_0, . . .)$. The intersection of a certain n-ary context is where a function of articulation has determined the partition of null in which certain individuals may or may not emerge.

To say x has the property of being a member of itself is to say of x that x has that property. However, to say this is to say that x has as one of its members something that is not-x and, therefore, x is such because it has the property that is not-x. As a member of itself, x is just in case there is that property defined by not-x, leading to a contradiction. In order for x to be the case is so just in case x is not-x, for x can only be such by that external property. The function of that property we have shown above leads to the identity paradox, concluding that one is x, because of x, as x is only had or understood by its method of construction. As such, it is related to itself by that property, a statement about x which is the context of x's assertion, A(x). A statement about that context cannot be proven on its own, so is added to it making $f(a_0)=A=$"a", the context (A) in which the name (a) applies $f(a)$ which populates the objects appropriate to each successive context of application. So, $A_0(a,x)=x$, $A'=A_1(x)$, and so on such that $A(a,x)=f(f(a),x)$. Thus, we model successive applications of that name. Each iteration a more complete statement of x than the last. Based on these operations, we can suggest that blackness' cultural endowment is embedded in each context of its assertion. Subject-continuity as a vector, lends itself to the ordinal construction of contexts. Each successor-context is valid based on a citation of the former. As such, blackness is no one thing but a movement. Blackness is blackness and not interchangeable with another thing.

So, the nominal-context A is the set whose only content is $\{f(a)\}=\{f(a,0)\}=\{$"a"$\}$, the assertion of $f_a$ at its zero. $\{$"a"$\}$ is the context in which that name applies, obtains. The initial assertion is the context of assertion. If a name is being used has a set of applications, we can assume that there is a function of use thereby acquire that function in a form of life zeroed for us yet indexed within our repertoire as that specific function discernible from others acquired at different times in different contexts. When applied, we get $\{f(a),\{$"a"$\}\}$ which is to say $A_0=$"a" to $A_0'=A_1$, etc. A's only content is a function, thus A(a,x), where x indexes the context of the assertion "a". So, A takes subject displayed by an operation of x and object a for its arguments, whose function expresses the identity of x as "a" in the context of A and its successors. The object "a" is the concept expressed by $f(f(a))$, the successor function of application. This model of a nominal theory is constructed based on the discussion of presets and prefunctions prior to application, existence ratified upon application, where the successor-context constructed by nominal assertion cites prior context of application to secure that the object of that assertion, appropriate to context, obtains.

A(a,x) such that we assume A("a=x") which contextually is just A(x,x). A(x) states the context in which x is a member. a=x in A is A given x is a. We put forward that the context in which an object exists is named, not the thing itself. These are the conditions in which a name's use is applicable. This is only with the possibility of being correct, for use is restricted to the chain of prior application which that assertion cites. Not necessarily about correctness as being correct is about the overall model qua reference frame. A name is meaningful because applicable within that framework qua context, even if the frame is incorrect. Get rid of inconsistent frame, the name no longer applies or should be used. The function of the application of a name we will detail below is as follows. I know that these conditions of which this name obtains the appropriate object. It is dubious to say "of" a thing "that" it is the thing satisfying the name.

Of/that is contingent, thus not compatible with identity. It is not necessary of/that, i.e. of me, that Victor. That/of necessarily holds within lines of context citation based on ordinal construction. In this way, we do not have an unrestricted extension of names such that one could say "of" a context "that" a name applies from a position external to that context, making that name applicable everywhere, therefore ambiguous. We say "that" this is the context "of" which the object named is a member. One cannot export a name into any context they want, amounting to an act of quantification into a domain without being a member of it, thus knowing what is in it without being in it—see Kripke (2011) and Quine (1956). As per ordinal logics, names are indexed to appropriate conditions via successor applications as a function of contexts.[9]

As an operation, subjectivity is revealed by some relation obtained between input and output that indicates some functional activity of subjectivity present but not represented in the set of inputs and outputs experienced within a particular state of affairs. However, a change in relations of output, soon input for this recursive-operation, expresses a change in subjectivity, indicating a subject at work. If we assume the basis for the function of identity indicates some subjective operation, we consider subjects in and of themselves as tautologies. This is to say that no one thing about x, $A(a,x)$, is such that of x, x is A, the context in which x is situated. Subjectivity embedded in null is expressed by the process in which identity emerges in the context formed from b to a above.

Let's say that there is an identity-function F that indicates a subject by relating input and output as the same individual. Let's call this the tautology, $F(x,x)$. Tautologies, then, are the basis to say any thing, thus to say some thing of x is to say that x is F related to itself. The content of the operation asserting "x" is the function animating what is known, called, F. For example, for all x, there is a property X such that, $X=f(x)$ holds true, for 0=0, 1=1, etc. for $x \rightarrow f(x)$ and F is the description of that function such that F-ness is the expression of the operation characterized $F(f(x),x)$.[10] In this way, F for all x states that $F(x,x)=F(X)$.

If input, from repertoire, and output, identity, are different, this operation entails a change in identity, the function of which entails a change in context,

---

[9] Is it necessary *of* Victor *that* he is Black? A statement equivalent to for-all/there-is quantification that may return a vacuous or invalid result. Or, is it the case *that* Victor *of* which is Black? Making the identity contingent upon definition, i.e. the reference frame constituting the context in which the function of that assertion returns the appropriate object of thought? The function of applying that identity has no means of justification save for being within a system of categorization that is not necessary to the entity that is labeled Victor, for the function of applying that label obtains only when utilized to construct that context. However, because "Black" has a function, its determinate before application allows us to understand that there is some operation, some subjectivity, that exhibits "Black"-ness reflected in B/~ intersection with other categories, thus partially grounding the expression of many. Because it utilizes the zero of that function, that subject operates outside of an identity with the label of a category and, yet, obtains within the form of life expressed by the affairs of subjects whose use of that function expresses various selves in the worlds in which they participate without having to predetermine the extension of B-ness. That individual takes up the functional content of B, utilizing it in ways that do not obtain given categorical stipulations while retaining that function in its endowment. Each act that does obtain a context of expression exhibits an aspect of B-ness unaccounted for within the identity-categorical schema that strictly defines the parameters of B. Those parameters themselves dependent upon the frame in which B is member. As part of multiple frames, different aspects of B arise, not all. The subject retains the zero of each function it acquires, utilizing those to which it was introduced to construct contexts or by being labeled within the context of others, thus co-constructing that context.

All have an initial state, a zero function qua mechanism, making subjectivity a means of expression by this operation, which is triggered once placed in particular conditions making for a particular subject. What is exhibited is based on application by subject itself or other subjects. Functions can be determined without the operation that runs over the set, whose objects are those functions, being finite thus its output predetermined. A capacity exhibiting infinite use of a finite being's means whose products are denumerable but it itself is nondeterminate for it is an operation of which the functions determined in context only define that subject after its application of those functions which is not a determination, definition of that subject qua capacity, but the function, the identity asserted, of which that same subject has many. The function is apparent only when applied in context but this does not imply that the function was not already present otherwise no thing would have been expressed. Without input no output. Initial operation merges function of self-assertion with functional content of context. From then on the configuration of functions evolves based on the form of life. Most shared for functions cannot stand on their conditions own thus the initial form of life creates what we considered humanity within which minor differences comprise unique individuals partially grounded by that initial state

[10] Kurt Gödel would state in his analysis of the logic of set-theory that all that is x, the set x, is only conceived by its method of constructing x expressing the property X by which x is related to itself. All else is a description of the property exhibited by x which is not a member of the set x, but how it is characterized, averting Russell's paradox. (see Gödel, 1940)

F(x,y)

That previous identity is F related to another in and only in this context. The failure of a particular identity's application within context can be considered

F(x,not-x) or F(x,y) and y=not-x

These represent contradictions of the identity paradox.

This system is unified by null, the set of no one thing and therefore everything which is considered the foundation for any articulation. We can model the construction of contexts with minimal stipulation by a negation operator indicating what lies external to the domain to which operations have access. In so doing, it is not necessary to define what is in the context, which would leave us with the problem of external stipulation once again. What is not necessary for consideration is bracketed and set aside, thereby framing the context required for consideration without having to rely on identities which are ambiguous to articulation at best, contradictory to the operation overall at worst. Taking a Negation-operator over null, it has been shown—by Geach, 1981; Hintikka, 1956; Peterson, 2018; Rodgers and Weihemier, 2012; and Wittgenstein, 1922—that any operation can be constructed by that N-operator.[11] This simplest operation is also clumsy for expository purposes, explain why it is overlooked.

If a system is built up by conditional assertion we can reduce the system of connectives with which assertions are constructed by a single operation, N. The function of this N-operator, here a joint negation, expresses the concept of null. In so doing,

N(x,not-x)

expresses the validity of the conditional "if x then x," it is not the case that x and not-x. We also formulate the articulation of possible alternative states of affairs. As null is part of every articulation, it is the basis from which that articulation was constructed and remains apart from each articulation as the set of any one axiom from which different models can be constructed. The individual x must exist in some sense if it is negated, just not in the parameters of the context articulated. Just as well, if null is at the core of each articulation, the joint set of all possible tautologies, it remains embedded in the chain of successor-descriptions built off of that core, regardless of its superficial appearance.[12]

With this understanding, we come to,

F(x,z), G(z,y)

so that F of x and G of y are related by a constant z.

---

[11] For example, an if/then condition in first-order logics is defined as the negation of a true antecedent and false consequent. Our successor-function above can be defined using this operator: 0, then if not-0 then 1, . . .; which is to say N(0N(N(0)N(1) . . .

[12] In this sense, blackness' subjectivity is an operation maintained across contexts. Its identities under the labels Black, et. al, are separate but not mutually exclusive from the function of subjectivity. Identities are context dependent in accordance with the ordinal construction of joint successive application modeled above. For within those contexts, to be "Black" entails an operation that pulls from a prior repertoire of identities. Thus, blackness is maintained despite its superficial appearance. To be "Black," or of any of its derivatives, entails so-called "mixed" antecedents. If, in institutionalized programs of Black study, it is maintained that the construct of Race is not interchangeable with a singular phenotype, then this should be uncontroversial. However, the controversy that surrounds this assertion indicates the intention to assert one thing, but hold on to an inconsistent position to maintain one's argument for future benefit. Otherwise, if the problem is dismantled without denying its existence or just revealed for what it is, its internal contradictions becoming an axiom in a subsequent framework for we know that argument is not the case, so would their position. The position would be on shaky ground for the status quo maintained would be shown inconsistent, self-defeating. This is not necessarily saying that the field would disappear, but would be better for it. We lose by picking which antecedents to call upon if and when convenient in order to concretize something, which by definition, is movement.

A system that determines z as the output of some function of x and y. The content of z must be *f*(z), its zero, Z, i.e. its Z-ness. The description of that system we would say is Z such that Z(F(x),G(y)). The same object obtained by the function of different inputs whose operation indicates the functional content of their assertion.

If F(x,y) and G(x,y) then F=G, so different operations can express the same function. The same subject qua capacity can be identified in multiple ways. F and G are merely descriptions of the contexts in which the operations of that mode of expression, subjectivity, runs. F(x,y) and F(z,w) are different operations, although nominally, they are the same function. As such, they indicate an effort of categorization.

F(G(xy)z) represents a concept we will term overdetermination. The function constructing that F-context determines the operation of another in order to set the possible roles obtainable by individuals within that state of affairs. A description of a description, F of G, is such that F does not have access to x,y but defines how G describes their relation, its output. The function of a function, in accordance with Frege, expresses a concept; it expresses an operation's role, its "function," in constructing a reference frame. This should not be too fast of a definition as we have shown how this is plausible based on the successor-function. We can also think of this based on ordinals. (Turing, 1938) If null unifies the zero of A, B, et. al, then each successor of A, B to some n-th place remains in that line, $A_1$. . ., $B_1$…, and so on. A function of that successor-operation causing one line to be a subset of another or for them to intersect, causes a pair-wise partition of those two lines in which an operation is only effective in and only in the scope of that newly constructed domain. E.g. F(x,x) implies F(X) namely that F labels the property X of x because x on its own is only conceived by *f*(x) and "x" of x is just x. The function of a function expresses this concept by producing an image. Here, if F of G, then F is the domain limiting the scope of the materials available to, therefore shaping the image of, G.  These expressions are then incorporated for use in the context of F. This is indicated by an operation of z indicating the function F. The terms of G are put to a use other than the one they acquired in their initial context.

How does this relate to the construction of a particular endowment brought with the assertion of any one identity? The immediate output of these operations populate our state. Changes in the relation one context has with others, that identities obtain with others, indicate the operation of some subject utilizing the means of their endowment to produce that identity in a way that is not necessarily caused by context, but, also, not random. It also accounts for why an identity may not have been experienced in that state-description. The reality qua relation between contexts of assertion was constructed by utilizing the means available in a way appropriate to the context of that assertion. Contexts and the state of affairs composed by the relations between them are constructs, and what is a member of one may not be of another. However, an identity can be introduced to others by a subject participating in more than one context. The functional role of that identity is dependent on the materials available in the context to which that subject is introduced. Thus, the same subject has multiple identities and its endowment exists across contextual boundaries. The endowment itself is composed of previous identities, functions now "features" of that endowment. It is by the function of these features that subjects gain access to other contexts as they express properties which are defined by a category which is context dependent. As such, previous identities can be applied in different ways based on the features of one's endowment. The subject utilizes its endowment's contents to construct contexts. The properties of others project a "self" in the shared contexts of that state described by their participation in co-constructing that context. It is well known in the philosophy of mathematics that every axiomatic theory, here subjectivity, admits of an unlimited number of concrete theorems, here identities, besides those from which the originating identities were derived, here an endowment. A cultural endowment as a form of life is therefore the domain in which this particular function obtains. The aggregate of these functions are the "features" of the endowment. If modes of expression are ways of being in the world (Baraka/Jones, 1963; Wittgenstein, 1953) then that form of life is indicated by a relation obtained between assertions as arguments of functions of construction.

If x is a particular set of features, "x" by itself is a tautology, it says no one thing on its own but potentially may be utilized to do so. We can assume if "x" is to be expressed, in the least it is an argument of its own function asserting itself: $F(x,x)$. From this it follows that every individual implies a mode of expression. This mode of expression is indicative of the form of life from which it was articulated and relates the arguments of various individuals which are distinct but not mutually exclusive. We have shown that modes of life can be considered predicates by their function over operations. $F(x,z)$, $G(y,z)$ for y implies $G(y,y)$ so a state F-G organized around z relating x and y. That state description is a relation between contexts composing a state of affairs which is only valid when the ordinal construction of the contexts that compose them are valid, the function relating contexts as arguments. As a constant, "z" represents how these individuals are related within the state of F-G.

As the function of the relation between operations as arguments, F-G are the apparent identities of the subject within a context whose parameters are understood by the limits to the function of F,G available in that environment. These operations obtain the same object, so we say that the function of contexts F and G are equivalent. Equivalent in the fact that they indicate similar operations but are not interchangeable. From this we can posit a function expressing the concept of this operation as $Z(F,G)$, a description of the operation of z expressing Z-ness. F-G relates the operations expressing x,y. It follows that if the set of the features of x,y obtain the same relation to z, then Z obtains the features of both F and G. We shall call this obtaining an affinity. We do not require an identity relation for it is sufficient to show the construction of the class Z and its members.

We will do this by saying that ~ represents a relation between arguments. In this way, "z" is the property of being labeled Z. So the relationship expressing the features qua function of Z will be illustrated as $Z/\sim$.[13] We should take care that we understand that the "features" of an endowment are functions, the objects of the endowment exhibited by the operation of those functions appear as the features whose relation within a particular context project the identity applicable within that context.

From here it suffices that we can prove a cultural endowment from the individual arguments that enter into such an operation. The function of these operations express subjectivity. The expression of this endowment is had by the operation of subjects that relate the objects of their environment in such a way that actively constructs the context in which that identity is asserted, placing them in the state of affairs in which their individual contexts are related. The relation between predicates attributed to this operation represent the relation between individuals, the structure of which is the state description. Thus, in accordance with Henry Louis Gates Jr, (1988) this endowment becomes a repertoire of signifying practices, utilizing the shared features of that context by some property that is not a member of Z and, thereby, is representative of another individual's repertoire. Those features not of Z and not of that individual are of some other individual that is a part of neither context. By virtue of this, the null set unifies the system.

We show this in the following way. If we have x such that $F(x,x)$, then the context expressed by the relation obtained by functions whose operation construct $F/\sim$, understood by an operation that gives us features whose function is f in x. Each f in x gives us a feature whose operation expresses an [f] of the context, form of life, of F. We represent this by $F/\sim[f]$, which is to say that there is a function whose underlying operation expresses F-ness. In this way f~f' can be shown for an individual x, the active construction of the features of x, which go on to give us all [f]-successors in F. Effectively, experiencing a change in the [f] of F resembles a change in the f's of x. The changes in input/output claims displays the function indicative of the operation expressing subjectivity without having to predetermine or stipulate its being interchangeable with something that it is not—see Appiah (1985). Changes in relations between the outputs of this operation, thereby, making a change in the context dictating which identities are applicable and which are not, model some change in the aspects of the subject constructing that context based on the inputs available. That subject is embedded in the contexts it creates. An "individual," then, is expressed by the relation holding the features of a context together, not the features in and of themselves.

---

[13] For a more detailed description of this concept, see the notion of presets and prefunctions detailed by work at *n*Lab (2019).

An identity operation that seeks to construct a system of categories in order to determine whether or not a particular aspect of an individual is a member of that state of affairs, thereby is that individual, fails. We know this from above, as "x"=X and F(x,x)=F(X) which cannot be proven on its own for by definition X=$f$(x) which makes for F(0), a label of no one thing. A relation in the cultural endowment does not tell us whether that identity holds in the state of affairs constructed. F/~ is incomplete as far as the decidability of categorization goes. The members of F/~ are present in Z/~ but are not necessarily represented there. We cannot determine, save by some other operation, if f=z. As a matter of fact, [f] cannot be identical with [z] for then we would run into the identity paradox. With each stipulation of a feature's being identical to an individual by the category denoting some contextually determined property affixed to that feature, there must be another feature generated that represents what the former is not in order to secure that determination. In turn, this occurrence makes every stipulation more ambiguous as more features are definitionally fixed for some property, determining the category of which that individual is a member. For each name affixed to feature, $(n^2-n)$ features must be proven outside that categorization in order to prove its validity.

Identity cannot prove itself on its own. In substantiating a singular identity for all, the required number of individuals outside of that category in order to the frame members of that category to make n=1, a universal stipulation, is $(1^2-1)=0$. For a single category given no one, $1/0=\infty$, and as such is non-determinate. So much for homogeneous identity claims. It is more likely that Z/~ would cannibalize itself as it constructs the context of its assertion based on what it is not. All the material supporting its assertion would be called into question or negated in order to secure an identity imposed from some other context external to the context currently inhabited. The relation exhibited by the operations of x can only be obtained within and only within F/~ as each [f] can gain access to F or F' by virtue of being able to show that f-successor constructed in x from its endowment. Since not-f and f can be proven within the system Z, making f and z interchangeable so as to complete a universal description of that system, this process leads to a contradiction in the overall onto-epistemology (da Silva, 2007) of that state of affairs. Overdetermination, then, cannot be complete due to the identity paradox. The mere instantiation of an identity whose dominance indicates its inability to be proven within the system of categorization, because it organizes that state, requires the production of that which can disprove that position lying just outside the framework imposed. For without that negated region, no framework would be discernible. Hence, our focus on the capacity to articulate forms of life, not the holding of some thing that is to be definitive, regardless of context, despite that form of life.

Take the collection of functions as objects of x's mode of expression. As there is no one feature interchangeable with x, then any feature of x can be conceived without having to quantify each one. Features are the functions qua constituents of "x" in whichever context, F, x is asserted. Call this feature $f_0$. We can describe a group by creating a set of those features without ordering them, call this {f}. From this we have $[f_1]$ of F. The next would be {f,{f}} in x which is $[f_2]$ of F, and so on. As there is no way, save by stipulation to know all the features that comprise x, we make it a nondeterminate, yet denumerable, collection. The total features of x is a known unknown, which is not to say that x is indeterminate or unknowable as it can be shown constructible. This collection of all features displays the concept that x's endowment is null for not a part, as of yet, of any one context. The endowment is composed of the zeroes of the functions qua features of and to which x has acquired. Thus, by negation we have an affirmation of this qualia without a stipulation of quantity.

Take the following.

$F_0=(f_1...f_n)$
$G_0=(e_1...e_n)$

$F_0=G_0$, just so long as

$f_0=e_0$

and, $f_1 \neq e_1$

Effectively, we have described an axiom of foundation. As there is a $f_0$ and $e_0$ as a member of each individual but not a part of their type distinctions, we can say that outside of identity there is a subject and hence a cultural endowment unified in null. That subject can be a part of F or G and yet remain apart from them in order to enter into alternative relations. This is the case for subjectivity is an operation, not a thing. Each feature gains that individual access to an alternative context which if taken from the contextual view, if that feature does not appear, we would say that that subject does not exist, their identity would be necessary to subjectivity. However, we know that this is not the case for x's role within different or the same state of affairs can be F or G. A single subject may have multiple identities dependent upon the context and yet remain the same subject. Also, that same identity, i.e. functional type, may obtain multiple subjects. It would be ridiculous, if not trivial to say that what is given is all there is. Everything may follow from something that is or is not the case; however, it is not always plausible that because something happened it must always be the case.

We can set up this simple operation producing features f of x's endowment by taking a collection of functions and calling that the set {f} in x of type-F represented, [f]. In this way, F has no scope on its own as all f in x are 0 in F for we have determined no one feature of F. By taking a set of that set, we have {f,{f}} of x in F which represents [f'] of F. The relation between $f_1$~[f'] in F is understood as a function representing the operation in x of {f} to {f,{f}}. We state that the functions in x, the ways in which x is F related to itself, are features of F whose overall function expresses F-ness by an operation in the endowment utilized by x.

We can do this for y, whose mode of expression through G is a feature [e] modeling e in y of type-G. In accordance with articulation theory, the aggregate of the pairings of features from x and y represents the intersection of their states of affairs at z. That state, call it Z, is the aggregate of those pairings, the structure of the context articulated from the set of conditions comprised by F,G.

As different operations may obtain the same state of affairs, articulation theory states that,

$Z/\sim[[f_1, e_1],[f_2, e_2],...]$

Therefore, rearticulation is such that G' is

$[e_0,...g[e_j...e_m]...e_n]$

Z becomes,

$Z/\sim[[f_1, e_2],[f_2, e_3],...]$

The operation indicating the subject is not only an extension of subjectivity identified within an appropriate context, because co-constructive of that context, but is also a deepening of the endowment. Within {0} from 0 we get, {0,{0}}, {0,{0}, {0,{0}}}, and so on, illustrating that within that initial finite context, there are plausibly infinitely many uses within that partition of null utilizing the functions of a subject's endowment. We also show that from 0 to {0}, {{0}}, {{{0}}}, {{{{0}}}}, and so on, demonstrates the extension of that subject's endowment across contexts. That progression is equivalent to 0, {0}=1, {0,1}=2, {0,1,2}=3, etc. such that we have $<n<x_1...x_n>^n>$ equivalent to $<<x_1...x_n>>$. Each x numeral represents a context citation, the members of each set the extension of the subject's endowment across and based on the prior construction of and participation in those contexts. Adding or taking away 0 from a context does not change the identity of that context as long as the particular framework remains intact. Thus, it is safe to assume that null is both a part of and apart from every articulation is uncontroversial. The operation of subjectivity is displayed by the active construction of these structures, not what is housed within them or by stipulating what those structures are which nevertheless relies on a frame of reference externally applied to the object of study.

As a result, we find that as Black identities become stricter and more narrow, blackness' subjectivity expands and deepens. Those depths are not altogether visible, accessible, to categories. The efforts to overdetermine and secure the referents of those categorical labels to ontological commitments are not definitive. Blackness grows exponentially with each attempt to restrict its expressions in different spheres. Much like Z-ness, we can show that B-ness is in excess of any finite categorization as it is exhibited by its method of construction and, dependent upon context, may appear under a different category housed within the endowed capacities present by not represented in any one of those finite categories. White homogeneity violates these rules of construction in the very act of its assertion—see appendix. Its extension is exhausted in one context. To extend a context over all others has no means of justification, for what is required is in a successor-context that, when evoked, automatically calls into question the validity of the domination of that former context and its ability to stand on its own.

Overdetermination cannot be complete for it is a concept relative to context. As a recursive operation, the articulation of subjectivity from a previous and continuously constructed endowment, based on the accumulation of prior functions developed as the capacities of that subject introduced, co-constructs the contexts of which the relation between subjects structuring our state of affairs. That endowment must be in excess of the context of any singular identity function's domain. Otherwise there would be no base for future articulations. We can think of the endowment as the residual from the identity paradox. The repertoire, then, is constructed with each instantiation of an identity. The context of its assertion can only be maintained if there are objects outside of the domain of that context, proving its efficacy as its domain is determined. In fact, we know that the repertoire grows at a rate of $(n^2-n)$ for each context determined by the assertion of some identity (n). Just as well the function of any proposition can be built up using only the N-operator, whose function is based off of null assumptions. If not, identities would be static and applicable everywhere, thus useless, not uniquely determining any one thing.

The operation of "x" indicates a particular endowment present but not represented within the particular domain or context in which it currently functions. We have also shown that the concept of overdetermination, predetermining the value and therefore the functional role of an individual from one context within the domain of another, cannot be complete. Through this, the failure of universal identity stipulations based on the assumption that an individual does not have a particular endowment is seriously flawed, if not self-defeating. This has major implications for the study of blackness, B-ness is not necessarily defined by B on its own. The operation expressing subjectivity is embedded in, because actively constructs, its cultural endowment.

As a result, the set of functions it contains, links assertion to context by the same input producing various successors as output. Each feature of the endowment indicates another property with which that individual is related to context. Thus, the same individual can become a member of multiple states of affairs without exhausting its subjectivity's repertoire. The use of endowed functions express a part of the stock of "self" via the identities available, yet portions of its self remain apart from the context presently inhabited. The self in itself is only had by an operation and therefore is this operation, within, because constructing, each successive articulation of its subjectivity. Without this, there is no hope of explaining the active construction of a form of life that does not presuppose the identities therein. If presupposed, those identities, without an explanation of how they came to be in the first place, are vacuous.

We say this in order to state that there is no need to overdetermine what B-ness, e.g. blackness, is by predetermining who does or does not fit into a category based on some property only known by the function that property obtains in expressing a form of life. The expression of that form of life is not necessarily what the label denotes. For those under that labe,l it is inevitable that they can be found utilizing or under others while remaining themselves. We need not rely on identity theory in order to formulate a theory of subjectivity. Stating that membership to category determines existence rests on shaky ground. No account is made of how the system got there in the first place in previous identity theories. Relying on such a theory presupposes a theory of subjectivity. As mode of expression is subjectivity, it is an operation which implies a method of constructing those expressions—infinite use of finite means. A method of constructing an explanation was shown in order to detail what F-ness is in the first place. F-ness cannot be "F," for F would be so just in case it is not-F. So each feature of x being put to use to express no one thing of x, but possibly some thing, automatically allows us to ascertain that feature as one of type-F. In fact, the operation expressing subjectivity populates the worlds, context, in which the identity as the function those objects of thought, exists. In showing how categories can only be grasped by their method of construction, which in turn implies an operation of construction and, thus, relating the entity constructing that system in the system that appears, we obtain a sufficient theory of subjectivity. As a bonus we get the underpinnings for an attempt at explaining how identity works without having to presuppose it by relying on mere stipulation with tenuous justification for assertion.

IV.

So far, we have shown how the internal operations of functions can be constructed. In so doing, functions become the "objects" of one's endowment, objects of thought. Intuitively, the move from no one thing to a set of things, and to a set of that set, resembles the move from 0 to 1, then 2. This resemblance indicates that there is some function present which is indicative of a relation between operations over two sets, one later considered input the other output. That function is not a thing present in the change. The function is revealed by the changes in input related to changes in output. An endowment is based on the construction of operations revealed as functions which later become predicates whose function are visible once applied. Predicates considered as functions—a description of x being F, such that x is F related to itself, e.g. $F(x,x)$—must have arguments. We can infer the presence of some operation that must be at work in order to construct the entity that houses these arguments. We understand $F(x,x)$ but if set to a, $a_x = f(a,x)$, thus the identity of a, is the relation between that label and the context of its assertion where numerals $x=1,2,3,...$ are understood by set, thereby context, theoretical construction—see von Neumann ordinals, 0 to $0\{0\}$, $\{0, \{0\}\}$, etc. which is $0=0$, $\{0\}=1$, $\{0,1\}=2$, etc. respectively. A set from everything possible can be a set from that set and so on, or another set distinct from the former, or a set of the set and those therein, contexts and types of contexts, or relations between them. If reference was based on the exchange of name for referent, $a=b$, the method of application vanishes and the identity or name is applicable everywhere and nowhere. We understand F within it context of assertion by its function f which indicates but is not interchangeable the with the expression of F-ness, the subject of that expression is the function of expression.

Functions as the content of one's endowment construct contexts of assertion. Predicates are the function of the arguments available within a particular context that individuals take up and put to use to project an identity that becomes an object in the context of others doing the same. This is how contexts come to be related in such a way to compose a state of affairs. An endowment composed of functions, *f*, as the content of its expressions is considered under the scheme,

*f*:=a→*f*(a)

so that if "a"=*f*(0), *f*(a)=a' when applied to some purpose. We state that if an object is plausible, there is a function of that object's construction, otherwise that object would be nothing. Changes in the relation between a set of inputs and changes in a set of outputs determine the role that object plays in the construction of our state of affairs qua each other's reference frames. The objects that make up that state description indicate the operation of a function, even if that object is a function of its self. A context is known by the projection of an individual's identity as a function of their endowment within the context of that identity's instantiation. As such, a diagonalization argument can be developed where as a table is constructed such that functions are subscripted to infinity, as the rows, and variables, as columns. This proof shows that the universal, self-referential, statement of systems cannot be proven in the systems they export. This also ensures that different cultural endowments built up by an operation over functions applied when applicable, because constructive of the contexts of their application, are possible when a diagonal is drawn across this table. With such a diagonal it become obvious that $f_0(a_0){\neq}f_1(a_1)$, i.e. $f(0){\neq}f(1)$. The unprovable statement is represented in each system indicated by the row. Each column is not empty but may or may not be put in relation with others so as to make a composite, co-constitutive, state of affairs. It is not necessarily the case that constructs are such that they overdetermine the output in one context for the sake of valuations in another.[14]

Implicit here is a critique of behaviorism's focus on output alone and its assumption or negation of this prior operation. Behaviorist disregard how output emerges in the first place. This course of actions misses the point because a function is no thing. Behaviorists' focus on external data does not account for how one organizes that data and assumes the system of categorization as implied in the data. Inevitably, this brings us back to the same problem of how those descriptions were constructed; back to the question of how and what produced the data externalized just as how the analyst externalized the system of categorization employed? On its own, a function says no one thing, but placed in the appropriate context because constructing that context, produces the correct output qua context thus conditions for appropriate subsequent input for reference to obtain in later assertions.[iii]

As modes of expression are the modalities in which beings are in the world, a cultural endowment, in accordance with what LeRoi Jones/Amiri Baraka and following Wittgenstein, is a collective improvisation. Improvisation is no longer an ambiguous term. It is the capacity to make infinite use of finite means, appropriate to context because constructive of contexts, therefore, not random but indicative of subjectivity. In this way, blackness is not interchangeable with one and only one thing within any historically constructed frame of reference. It is also contemporary, immediate, as blackness is the capacity to express many selves yet retain its subject's continuity.

This comes with the idea that a form of life is such because of an affinity of signifying practices that lies outside of the intervals categorized over the arc of that mode of expression. The identity paradox comes about when an interval is supposed interchangeable with that arc. One must participate in the active construction of that form of life. If not, then to quantify into that opaque endowment from an external stance misses the mark of what is actually there. (Quine, 1956) Blackness is an affinity between many modes of being, which if we do not rely on identity, other forms of life may enter into an affinity with it. However, to quantify into blackness so as to determine its value for use in a context in which it is barred from actively constructing does not make the one doing the quantification "Black." That identity is not the operation expressive of blackness, but makes that entity negate all that makes that identity what it is. This reveals an attempt at dispossession so as to construct its own context which will remain void, for that identifier is at the expense of what animated it in the first place. As contexts are built ordinally, this act cannot be disavowed, accountability can be tracked by the only content of the propositions of that external position. Its function is constructive of that void because it displays an extractive system of categorization, hence control.

---

[14] Diagonalization states that for all n labeled N, if there is a function $f_n(n)$, it would follow that $f(n)=f_n(n)+1=f(n)+1$, a contradiction. Gödel showed that this can be conceived because it can be constructed, computable, by coding the inconsistency. By Turing, that inconsistency would be stored as an axiom from which a successor frame of reference can handle this concept because it was used to construct the frame in which that function is applicable. Gödel was concerned with encoding self-reference, hence our theory of the assertion of names as only having functional content. Our theory takes the form $f(a,x)$="a" so that x=1,2,3 which are the contexts of application based on ordinal construction. After Gödel, we would say if "a"=A, then in A, we have the function for which "a" applies and the definite operation paired with that function of application when appropriate indexed by that method of construction, i.e. computation.

Why do certain features whose function express a form of life come together? From external stipulations. However, although grouped in such a way, because of the subjects capacity make infinite use of finite means, a form of life develops which is not at the expense of the individuals therein. Their modes of expression are different and remain such, but reach an affinity which co-construct that form of life. They relate the contexts of their assertions in such a way that constructs a state of affairs that is not a monolith, indicating a people that is not such either. In so doing, there is the capacity to rearticulate states of affairs without denying the historical conditions of their construction as well as without denying the creative potential in the multiple realities represented. These are real because modes of expression harbor consequences for the one's doing the construction. To actualize different or at least one of those realities implies articulation. "Actual" is its function within a reference frame. The issue of how reality and actual interface is one for contemporary study, however, in rearticulating that frame, since looking through a frame one cannot picture the frame itself, we change the framework in which actualities emerge. Correlation has little to do with causality, identity is contingent upon subjectivity. Identity did not necessarily cause subjectivity however that subject can utilize the finite means once placed in that group to articulate a self not necessarily caused by their historical condition but without denying them; thus rearticulating that groups relation to and with the others that make up that system. Those defined within that system can temporarily benefit from it. As systems are constructed the work is done towards the benefit of the one who constructed it. Since the propositions expressed from endowments only content are functions, they express nothing until applied and once applied relate the function of assertion to the context of application.[15]

Our theory is not the deployment of tokens or representations, the superficial appearance of what comes out one's mouth, but the deployment of functions which in each successive application gain their character. In this way, someone can utilize a phrase taken to be one way but "meant" in another. The use still is valid to the utterer, as the function of that phrase as its been constructed over successive deployments maps the way it is used to the context of appropriate application. Just as well, those others may acquire what may by its use further the strictness of their determined position in that state to the benefit of another no longer in the picture. Structured thus, when contexts of assertion are then related so as to provide the architecture of our state of affairs, the utterer can be held accountable for the function of that phrase as others acquire that function in their repertoire. To acquire that function is to have that concept such that the concept is expressed by the function of the function applied, which is to say that one has the function and because co-constructing the context, can also apply it. They have the determinate of the function F(x,x) and their assertion of *f*(x) is the input to an operation constructing a context related to others whose object is a constant within the state of affairs of the individuals functional construct of the contexts that make up that state of affairs.

---

[15] For example, an endowment is the software and as such is a bundle of instructions without input. Without input, no output; however, because developed, its uses become a repertoire constructed by the contexts in which and to which it is applied. From use, software becomes an application. If not applied, its assertions on their own are empty, therefore are not assertions at all. From an initial state, software is developed by an assertion which sets the first limits for the type of input sought after for future output, thus giving us unique modes of expression both from the endowment being developed and in contrast to others.

In sum, the smallest common denominator of all subject's projection of an identity in the appropriate context is 0, for that identity is only such by way of its method of construction from null. Thus, as that denominator, subjectivity is a generative faculty, an inner creative capacity, for 1/0=∞. Just the same, the common identity to all is null. Making an analysis with identity as primary is a wrong start because it is not a necessary constituent of subjectivity. Therefore, inter/infra-subjectivity is theoretically plausible because constructible outside of strict adherence to finite categorization stipulated as universal over subjects across contexts. Subjectivity is the capacity to make infinite use of finite means. If universally stipulated, this then is a reflection of what is required for the one making the imposition, thus a part of their identity, not the one identified. We arrive at a constructive idea of the concept of states of affairs and affinity. A state of affairs is expressed as a function of the relation obtained between the contexts of identity assertion, a coalition between identities. An affinity is inter/infra-subjective, a pair-wise operation between the functions as content of the expression of subjectivities utilizing an endowment qua repertoire of signifying practices. If someone attempts to quantify from one context into the endowment constructing another, as per Quine, that quantification either leads to inconsistencies or vacuous propositions. By showing how different endowments construct their respective frames of reference as well as enter into relations with each other, we show that it is possible that each endowment can maintain its integrity without incorporating content by dispossessing others of that which their form of life constructs, i.e. without what is most often nonchalantly proposed as and uncritically interrogated under the label of appropriation.

V.

The purpose of this essay was to attempt to develop, as best we can, a working theory of subjectivity in order to account for blackness as a mode of expression therefore its way of being in the world. In so doing, we would be able to explore and provide proof of the concept of a cultural endowment, the active construction of the subjective continuity of a people by virtue of the activity of its subjects even if from a null state. As a bonus, we showed that null is not void. Rather, it is the set of no one thing and therefore the set that provides the material to articulate any thing. As such, the possibility of constructing contexts in which identity can be asserted with the least amount of assumptions was determined by some operation over null. This was done to disprove a prevailing notion in the study of Black-ness within Critical Race theory called Afro-Pessimism (Patterson, 1982; Sexton, 2016; Wilderson III, 2008; et. al) This theory supposes that blackness is absolute nothing and operates out of a cultural death. The Pessimist is wrong for their determination of blackness presupposes its being interchangeable with a finite set of identities characterized by their axiological determination as absolutely nothing. Designating that blackness is only a predicate interchangeable with subjectivity fails when the logical structure of the Pessimist's claim is considered critically. So constructed, the Pessimist's determination is vacuous, x=F(…); or indeterminate, either x or F; or self-defeating, F(x) assumes x is, therefore not absolutely nothing. It is contradictory because blackness is expressed by an operation, it is not a thing that can be exchanged for a predicate that inevitably is void because devoid of argument. Therefore, blackness is a member of the multiple contexts which it constructs under different names. Blackness cannot be a void for it actively constructs a repertoire from the contexts in which it participates. It cannot be nothing for it is no one thing but operation that can go on to express any thing appropriate to context because it makes use of the means therein. As such, the Pessimist position should be given up. Even the use of silence is proof of the functional content of assertions. We know that F on its own is equivalent to *f*(0), the basis for assertion at all. Its use can show inconsistencies in the state constructed or realize different relations latent therein. The potential to construct alternative, correct or discard inconsistent models/contexts, lies at *f*(0), which is where contexts connect forming our states of affairs. We know this because current realities, and possible new ones, meet at their disjoint union in null, connecting us with others. The difference between mere static and music are the silences which hold the relation between utterances, revealing the syncopated output as song.

If we were to say, "for all Blacks there exists a label," this quantification is equivalent to, "there does not exist a Black for which it is not-not for all labeled not the case." The quantifier "all" is equivalent to "there-does-not-exist-an-entity-such-that-it-is-not-the-case," and "there-exists" is equivalent to "not for all not." This statement proves that there does exist some Black that is not under that label. Most likely, there are no Blacks that obtain that label for the formula returns not-Black Blacks or a label that labels no-(one)-thing. This Russellian paradox allows us to see how one can state something without that statement capturing any one thing. Revealed as merely a contextually dependent description, we show the Pessimist can talk about Black/African-American identity as being nothing without there actually being individuals whose operation expresses the form of life of blackness in such a way that obtain that description. As a result, we can explain why the surface appearance of the statement seems true, but is actually empty. The predicament of the Pessimist stands as the sign of an erroneous judgment bent towards maintaining a problem for benefit, rather than analyzing it. By treating identities as names, the same subject can have multiple names but in context they are the "subject" of only one. Contextual definite descriptions of individual show that a name is a label applied to a property indicated by a single feature of the many that individual possesses. Names then are not the individual but their use indicate the system of identification organizing that state of affairs. As an individual has multiple features, each indicating a different property that gains it access to multiple categories, identification is internally necessary within context, but identity itself is contingent to context when it purports it is not.

We have demonstrated this on two fronts. A description or identifier cannot stand on its own. This would be akin to a sentence without object, i.e. its contexts, and the object the function therein being the "subject" of this activity. The function of a description with absolutely no argument returns nothing. As such, it would not be a valid description and is self-defeating. Identifiers, names, rigid designations, etc. have no scope unless applied to, i.e. become a function of, some object. Thus, $F(x)$ is just $F(x,x)$ for x is F related to itself in the context indicated by the function of F. x is x and no other thing. Blackness is blackness and no other thing. A movement is movement, and not another thing. Hence our treatment of blackness' subjectivity as a vector so as to emphasize that point. $F(x,0)$ affirms x and its relative description as applicable in no one context. $F_0$ is a member of the set of all possible state descriptions. $F(0)$ returns 0 when F attempts to stand on its own—indiscernible, from an external stance, because the zeroed determinate as functional content of any predicate qua label, G,Z, etc., remains so to others until applied. However, the existence of these functions are still representable before application. Therefore, 0 affirms that that which is labeled nothing must exist, in some sense in order to be labeled at all. However, the function of that identity says nothing for no one operation is utilizing F to say something of its self/world. "Absolute nothing" is merely a description of a thing even if it disavows that thing in order to assert itself. What it disavows must exist in some sense, namely, not in the context in which F exists. Afro-Pessimism, then, is the vacuous proposition, not blackness; for blackness is an expression of subjectivity.

So, we state that $F(0)$ connotes the property of nothing which those that fall under some category after the application of the name of that property express the concept nothing. As such, nothing would be universally applicable, therefore useless. As null is the only property which has itself as sole member, null is both part of every category constructed, for they are successors from null, and yet remains apart from each category. What is in excess of the partition of null that falls under that label becomes the material for further articulations. If blackness is null according to the Pessimist, then everything follows from blackness. If an operation over null indicates the presence of subjectivity, then affinities between different subjects is possible. Divisions by virtue of one's discernment between identities are a matter of stipulation. Their function to some end, which may or may not be useful, but is not necessary to a working theory of subjectivity, is applicable only to the context under consideration. We also see how if starting from identity, alternatives can be constructed through the very same subject-operation, utilizing those means and displaying a transfinite continuity despite the imposition of finite contexts organized by the system of categorization in which those identities are available. This is so we do not deny the world order, the racialized injunction we have found or do find ourselves in currently. Null is a member of every set and yet the only set that can stand on its own for $F(0)=0$ no matter what the operation or the category applied. Thus, how terms operate is a more fruitful conversation than how they are defined. Definitions are dependent upon and only useful in and only in the contexts which they actively construct. As Amiri Baraka would state in 1965, Hell is our definitions.

The implication of this finding is the conclusion that the rearrangement of the features which indicate subjects being members of some endowment changes the property expressed by an operation over that endowment. Dependent upon the context in which the subject projects itself in that context, subjects become the object of the function of another's repertoire. We can make clear the distinction between identity and identification via the overdetermination of one's modes of expression, thus way of being in the world. Identification is a function of a predicate applied to a projection of one's subjectivity. As such, and in defining what that subject is within other contexts, the relation between the contexts of assertion composing our state of affairs, a definition of a definition over a subject, is not complete. The subject is in excess of any one identity because it participates in a multitude of contexts. However, once the function of these terms are the focus, their use towards some end become explicit so that terms previously innocuous are revealed to maintain and make more strict the system of categorization and the commitments it entails. Explicit now are the acts ensuring benefits of use for one context from the means generating whatever value to be derived in another. As shown above, the "content" of an expression, if any, is its function, not some thing, as indicated by the appropriate context in which it is applied.[16] The context in which it is applied indicates the concept which the behavior of that expression displays. Thus, an expression's use indicates the context of that use, determines that context. A phrase used in a "racist" way indicates a racist context. One organized by the commitments the deployment of the concept of racism by way of taking advantage of that system, the function of the identifier in that phrase, entails. This is so regardless of the superficial appearance of that phrase, as it is that phrase's use that constructs that context which appropriately characterizes its content. The function being the object asserted within the context indicated by the use of that phrase, under the label which it actually represents. That racist function is the content of that phrase. Change the property alone, one just reduplicates that endowment and its individuals under a different name but continuity is maintained.

In summary, an endowment is a repertoire of functions constructed and accumulated as one instantiates contexts by virtue of the assertion of an identity by some operation. Subjectivity is this operation that runs over that endowment of functions in order to express a self within the context constructed. Identity is the projection, output, expressed and the input for future articulation, here rearticulations. Subjectivity's assertion of a self actively constructs the contexts in which the assertion of that identity is appropriate. Each identity indicates the context of its assertion, the relation between these contexts structure the state of affairs of the individuals who actively participate determining that state.

This leads us to the importance of this finding. To assume that one can change the world before they change themselves assumes that one is separate from the world. This is not to say to not help people, but in changing world with the idea that one is separate from it, only treats the effects/products of systems without attending to the conditions in which those issues recursively emerge and thus doomed to repeat those same issues. The articulation of self, implicates one in the articulation of the contexts in which they participate. Changing the world before a change in self reduplicates issues of self in the world. That world is made in the image of, as a function of, the author(ity) assumed separate and independent from it. We have shown that this assumption leaves the assertions of that sole individual either vacuous or null. If null, assertion re-implicates them with the others participating in the construction of that context. The transcendentalist become the means of their own subordination. One is themselves because they have the property of being in the context which they actively construct. One cannot disassociate themselves from the contexts they actively participate in and therefore construct, so it follows that trying to change an individual other than self in context, does not change context. Thus, one is only a self because in the world of others and for each other they are because of the same operation. This is what LeRoi Jones/Amiri Baraka called a collective improvisation. Others are because they're in my world; and I am because I am in others. A change in one's self-operation marks a change in the world one actively constructs by being in it and has the added benefit of changing the structures through which a world, our world, is expressed. Changing self so that context changes, alternative contexts arise and inconsistent ones return empty results, capturing nothing. In this way, one reaches affinities with others rather than subjecting others under a singular domain. This can work in the opposite way, however, because those working towards homogeneous structurings of their world, enacting a system for self at the expense of those whose operation provides the value of and for that world, those models fail to be consistent nor are they sustainable.

---

[16] See Appiah (1985) on assertion and conditionals alongside our treatment of the ordinal construction of contexts appropriate for the application of an expression as a function from the stock of the subject's endowment. In fact, Appiah develops his theory based on unified partitions within null. For example, we conditionalize on the extent that R is the case within the partition [(R),(not-R)] deriving all possible from the calculation p(R) such that R-[1-(not-R)]; which if, in his terms conditionalization, in ours articulation, the construction of contexts based on our N-operation utilized by subjectivity, is if R then R states not-R or R which is to say not-R+R, the disjoint union of if R is the case, unified in null.

As contexts are constructed by the assertion of "self," implying the world inhabited by that self, changes in context, marks changes in "mind," thus a change in worlds. Certain functions, identities, become inapplicable and ridiculous to assert. This is how we could attend to issues of racism and other forms of subordination. Alternative structures become necessarily possible. It is almost too much to ask to change minds directly. One cannot peer in, see and rearrange the stock of metaphysical things that constitute another's mind—see Kripke's Wittgenstein (1982). However, working on the self as a member of the contexts in which they participate with others, otherwise there is no "self" to be had, changes the context and what is expressed. In this way, we have a better understanding of why that same Baraka would say that we must attack with a hundred shades of Black, not just one regardless of which side of the color line an imposition of a singular identity comes from. Based on affinity rather than coalitions, we state that affinity is not based on interests. Interests imply an asymmetry in structure, the imposition of dominant and subordinate positions. Affinity is a relation engendered by the coming together of conditions as the context of individual's assertions of self, independent but not mutually exclusive. Derrick Bell Jr. (1980) would analyze the convergence of interests as having the adverse effect of concessions made if and only if benefits can be controlled by the one making them.

Above we stated that an endowment is constituted by the zero of multiple functions, determinate to the possessor, but as yet experienced by others for they have not been applied. Starts from an initial state, a function that merges self with the assertion, the construction, of their initial world which is not-self; the disjoint union of which displays the functional relation of one's-self in that world. From then on successive applications of that function relates that subject with others, thereby acquiring different functions as objects whose application constitutes that subject's state of affairs with other subjects, contexts, and reference frames. This increases the objects, concepts, of its repertoire, that subject's endowment.

Implications of this theory to the concepts of mind, memory, choice, and knowledge are profound in this respect. Maybe it is the case that one can hold different, possibly contradictory thoughts, in one's head and apply, choose, neither? Maybe it is this that marks the uniqueness of the concept of mind, marks our humanity? It's not necessarily that we memorize "things" that reside in the extra-mental world. A finite mind cannot possibly retain a potentially infinite amount of images qua information. It is that we acquire and determine functions which become available when the conditions in which their application is appropriate are available. Conditions, the environments we actively participate in because we construct them, prompt an order of operations utilizing these functions. The zero of these functions are discernible to the possessor but not to the context until applied. Their application structure the contexts in which we can describe a person's having "knowledge." These functions are evoked in these situations because that subject previously utilized that function in order to participate in the construction of that context, thus determining that function, indexing function to context.[17] This operation expresses the concept of memory. These are the contexts we go on to say that those concepts, the functional content of an assertion, obtain. The zero of a function is the function itself, which is no one thing until applied so as to construct the context of its assertion, the object of which is the concept, the functional content, of that assertion. This makes sense as a function possessed without input generates no output perceivable by others. Yet, that function is still well-defined for the possessor and as such is discernible from other functions in the inventory, whose relation expresses the concept mind, here subjectivity, of and for the one who constructed and utilizes that function qua object of their endowment. Knowledge and reference, then, are artifacts of how we used to conceive this process—see Papineau (2019). A single entity can have multiple identities; a single identity can apply to multiple entities. Reference, then, is interminably ambiguous and the notion of what counts as knowledge remains foggy at best when based on the former notion. Functions are determinable, i.e. computable, while reference is not. (see Gödel's diagonalization proofs, 1931) A finite vessel does not hold infinite things, but with determinate functions, whose use construct contexts in which the objects of thought, concepts, obtain, alongside the zero of functions representing the function itself, we can make sense of this faculty of subjectivity as a capacity such that the concept of mind is displayed by its making infinite use of finite means.

In this way, we share the practices we develop, as states of affairs are co-constructed at the intersection of individual context's and whose relation express the world we co-inhabit. We also share the knowledge of functions that do not work. As such, all are again accountable for the state of affairs they construct and project; averting the predetermination of a state in which the subjectivity of persons cannot emerge for one has to trade themselves for that state, a future prefigured from a past determination and thus doomed to relive a recurring historical construct. A homogeneous state and a single individual over others are equivalent propositions in this respect. Subjectivity is ascertained by an operation indicated by the function of an identity or name, (n) whose expression indicates a context or world, (w) and whose function relates that name within that world, r(n), becoming that world's content.[iv] Here, in vector notation, for subjectivity is an operation, not a thing, and therefore only captured in process, in movement, the concept of subjectivity is displayed by the ordered-triplet, $<w,n,r(n)>=<w,<n,r(n)>>$, (Gödel, 1940) the world's content is the function of the names utilized to construct it.

There are many branches within the cultural endowment whose utilization articulates blackness. If a singular identity conflated with ethnicity despite the capacity expressing the form of life of a people presides over everyone therein, a notion of purity is implied. This purging of that categorical distinction does not arise when we come through the use of identifiers by subjectivity, for the same capacity inherent, constitutive, of all individuals makes it so that they're not mutually exclusive. The form of life obtained is recursively denumerable, composed of many whose mode of expression reach an affinity with each other. The objects produced do not necessarily have to be the same. Within blackness there are many "ethnicities." Blackness is not confined to one nationality or geography. Blackness is not interchangeable with any one, for then it would not be blackness.

If blackness is no one thing, then it is that which places systems of categorization into question. As such, blackness becomes a perfect candidate for analyzing how categorization and identity come to be by way of the various methods of their construction. It, then, is possible to make an account for subjectivity, which must come prior to identity and the categories to which they are tied, outside of the confines of identity concomitant with categorical stipulations. Without the existence of this prior operation qua subject, this inner creative capacity, there would be no such thing as identity, as construct or applied because constructed by someone else. If subjectivity is in excess of what it produces, then we can also account for the emergence of "cultures" and forms of life that function without a reliance on identity and persist although unidentified. Subjectivity is characterized by a capacity to make infinite use of finite means; the expression of a faculty that produces a digital infinity. As such, this operation is unique.

---

[17] Situated knowledges qua Standpoint epistemologies have been developed and are key to black feminist scholarship. In this way, we can understand epistemologies and their ontological commitments with regard to their ethical deployment without hopeless relativism.

Black culture being the primary export for some and there being contexts in which Blacks construct the only value therein, the logic of the Pessimist's claim has been shown false in addition to the empirical evidence to the contrary. Black is only an identity with an ontological commitment to there existing Black individuals. If that identity is cultural death then I say you can have that identity. If it is only true as a historical construct, fair, for its formation as such harbors consequences for us today, therefore is a reality. However, this position is still revealed as merely as a name in some other system. Blackness has been shown to operate outside of identity stipulation for it is a member of a multitude of contexts, that stipulation only being one. Its subject sustains itself and lies everywhere for it cannot be placed no one where. Naysayers should check the form of life to which they are affiliated, rather than condemn those they cannot see because they are no longer a part of it. This, especially so, if they hold themselves, according to Frank Wilderson, "transcendent" to that form of life. As transcendent, for assertions are only applicable in the contexts they actively construct, their argument is vacuous for the speaker knows not of what they speak. It is impossible for there to be absolutely nothing. Otherwise there would be no basis for what the Pessimist asserts.

Conflating an operation for any thing was the wrong start. The question in regards to Black subjectivity goes beyond retroactive descriptions of its value, whether those valuations are good enough or proof of how bad that marker is. The Pessimist course seems to have been drawn up in order to substantiate an already held belief. This belief only being attributed to the products of this capacity in the context befitting a claim made of one context and over an activity which has the capacity to produce both "good" and "bad" but nevertheless cannot be nothing. The very fact that the identifier Black has a function in our state of affairs indicates there is something animating that label. Efforts negating what animates that identity are dubious because they are self-defeating. The function of these terms indicates that there is something using them or something external to the context of their use imposing them. This automatically entails that a single label does not exhaust what animates that label in the first place. There are ways to build forms of life outside of the categories determined by these identities and assumed uncritically as the basis for whatever become the "cultural objects, products, or output of that form of life. The fact that we can ask questions about the frameworks imposed on these activities means that alternatives are possible, for if there were none, we would not be able to question them at all. If we could not make alternative affinities with one another within the forms of life we inhabit and in excess of the identities and their contexts we assume, we would not be able to form questions in regards to these supposed categories necessary to our being, our membership, in the world. We have shown by some operation how blackness from no one thing, or a stipulation of a finite context of finite means, may make infinite use of those means which inevitably constructs an alternative state of affairs and challenging that previous system of value.

[^i] APPENDIX

Operations within workspaces operate under the concept of internal merge within Chomsky's outline of the Minimalist Program for modern linguistics. What appears as "copies" for semantic interpretation in the deep-structure of what is externalized are the functions of the term initially merged with others. Our model: a,b; {a,b};{{a},{a',b}}; and {$f$(a),{a',b}}. The function $f$ is present but not represented in the successor-context. As the zero of a function, this can be assumed where $f$(a)=0 and by set-theoretical construction, ordinals are constituted in the following way: $f$(0)=constant; {0}="1", context of assertion; and {0,{0}}={0,{1}}="2," defining a relation between 0 and 1 in the context labeled "2," i.e. displaying the structure of function $f$ and making the function whose presence was there before construction discernable in the current context. That context is constructed by an $f$ utilized from the endowment of the subject constructing the context in which the application of "a" is appropriate, thus returning the value qua object of thought, b, that the function of "a" expresses. With the initial assertion of a's name being the construction of the context {a}, appropriate application concerning the context in which that assertion obtains is modeled, {{a},{a',b}} to {{a},{a',b},{a''{a',b}}}, and so on. According to the set-theoretical construction of functions, $f$(a) states that $f$ assigns a to the object b. If "a" obtains, $f$(a)'s use, in a particular context, the domain, {a}, in which its assertion applies, that nominal-function's use is appropriate because it constructs that context by that assertion, returning some value in the context in which the object of thought it names, b, obtains. The content of b as the assignment of the successor-context is $f$(a). That determinate as the function of context and assertion is as close we come to the referent of names. Names refer to the context in which their assertion is applicable, because appropriate, due to that function's prior and active construction of that initial to successor-contexts. The function remains embedded. The highest structurally, i.e. the most recent assertion, then is output in what that program calls the surface structure of the phrase.

We make sense of this as follows. If "a"=$f$(0)=name/constant on its own, then "a"={a,{a}}=f(a) the content of the tag "a"={$f$(a)} it's domain of application. So, {$f$(a),{$f$(a),b}}={$f$(a),{a',b}} where b=a' and the determination of $f$(a) given an operation utilizing $f$(a) is that function's application implying a context of application. The concept, the function of the "function" as the operation expressing b, is determined by the function of "a". We assume "b"={$f$(b)}, but in this context b is not the subject exhibited by the use of that function. It is the object of the proposition constructed by the function of the name "a." This shows "a" as a tag for the contexts that the mental object b obtains, which is just the successor context in which "a" is appropriate. If "a"=$f$(0)=A then {A{a',b}} which says $f$($f$(a))=b, so b=a', the context name in which a's function can be discerned. What is stated is "a," the successor-context b is expressed by the function of "a" which does not appear in the initial context because b is not physical. The content of the assertion in b is just the function of "a" despite its superficial appearance as b; the subject, the function, of that assertion nevertheless remains, expresses, "a" as the structure of the nominal act is a well-defined function of the predicament in which the operation of applying that name obtains. We can continue this operation via {a,{a,b},{a,{b,c}}}, so that the form {a,{a,x'}}, were x=successor-context, retains the consistent use, function as application of the constant, "a" based on von Neumann ordinal construction. Same name, function, the same object, argument, for the contexts in which that name applies. Same name, different objects, successor-argument, that fall under the category expressed by the function of that name's application to the objects, function qua content of context of assertion, having the property denotes by that name in the contexts in which that name use constructs and therefore applies. So, b,c,d,... are the objects of thought that fall under the category expressed by the function of that nominal constant "a" as those successors b,c, etc. vary in respect to the context that names assertion constructs.

Where b=a' just makes b part of the domain, the range of determinates of "a" is in accordance with internal merge rules. Above, when we showed a,b becoming the set, object of thought whose internal structure is {a,b}, we modeled access to one's endowment as access to the zero or function of that name. We represented, {$f$(a),{a,b}}, as {$f$(a),c} for {ab}=c so that "if (a×b) then c," where (a×b) implies $f$(a,b), then c enters an operation with the function "a", not a in and of itself. So if $f$(a,b) then c, then {a,b}=c becomes the object of f, in the domains in which "a" is applicable, constructible; $f$(a,{a,b})) is just $f$(c,"a") with the function of the name being used being well-defined for $f$:={a{a}} the one to one mapping of a to itself which does not violate internal merge. In this way the name "a" references the object b which is a stand-in for the internal structure of use {a{a}}, the function $f$(a) as the relation "a" obtains within the context of a it constructs. An expression of the applicable domain which is just the successor-context that is a built off of the initial use of that name is the function of that originating context of assertion, determining that function in the first place as a relation to a. The form, structure, of reference remains intact with each successive use of that name, regardless of the symbols we used to "formalize" for explanatory purposes: $\{a,\{a\}\}=\{f(a)\}=$"a"; where our notation is just a handy way of describing the relationships we wish to precisify; and as the determinate of a function is its zero, "a"$=f(0)=0$ so it's use is $\{f(a)\}$ and so $\{f(a),\{a,b\}\}=\{0, \{0,1\}\}=\{0,1\}$ which sets the function of names use as the context of its assertion.

Thus, in accordance with internal merge, the items entering into the operation are not encoded with information, requiring stipulations which cannot be justified. However, their content is the capacity to encode, construct, the contexts which we inhabit because those contexts are such due to the functions acquired by this capacity of that subject.

[ii] APPENDIX II

"If c then c" means that an operation over the domain in which c emerges is coherent in the context of its assertion; its mode of expression being the modality in which "c" is in that world.

"If c then, if (a and b), then c," is licensed because it is the case that c given (a and b); and since the operation (a and b) expresses the conditions in which c obtains, and "if c then c," then we are ok.

Take a and b such that they are comprised of elements,

$a_1,a_2,...a_n$
$b_1,b_2,...b_n$

There must also be an $a_0$ and $b_0$ that say no one thing about either and as such are the labels applied to each A,B. As a constant, $a(0)=A$, function without arguments, cannot be proven in that context because it is the context; what's considered a member of that context is defined on the basis of that identity, thus constructs the context from that null starting point. As every context has this zero, the disjoint union of them all entails that contexts are constructed from and within a non-empty null. Every statement made about the context of which "a" is part is provable by some feature within the domain of a. However, statements in that system are not provable in "a," for a is expressed by the function of that system's construction. A universally reflexive statement, the name $a_0$, is that unprovable statement included in a. The zero of every possible individual in the context of their assertion, it is the label of that context. That label by the function of its application creates the successor-context a' whose validity can be traced, references A, the reference frame constituted by the initial assertion of "a." Valid in its assertion, the context nevertheless houses what can be seen as contradictory. This is "a" but not-a, "not a" in and of itself but a function of/in "a" as use constructs the context in which that subject participates as the name applied to whatever in that context is "a." As a,b are our initial conditions, our null assumption is that a,b are unified by what is not the case for either: a, b, not-a, not-b are all possible in null as the not-operator reveals the method constructing both a or b from all that could be the case. As such, those labels represent our understanding that there is a function constructing each collection. A label is only such when applied; so for a, there is a label A for all that is a, such that $f(a_0)=A$, as the label of the category to which entities have access by some property determined by their possession of some feature. As all that is the case intersects at null, the property relation is the zero of the function by which a is constructed. Thus, that label says no one thing until applied to the set constructed by some function of application to that entity in the context of A. As a rule, that label's function is indexed to its context of application. We define an operation conditionalizing, i.e. constructing, from null what is the case by the relation, "less than or equal to"; an operation pairing features of those strings to articulate some particular out of all possible via multiplication.

The operation, "if a then b," becomes a≤b. This results in the construction of the domain $b^a$ so that $a_1=b_1$, $a_2=b_2$...We obtain $b^a$ which is to say a universal statement of everything extends as far as no one thing, for b is the case just so long as a and if not-b, then not-a. However, we have properly defined the conditions in which any thing may emerge. When $f(a_0)=$"a" the assertion states no one thing/feature/aspect in particular about a, only everything known/about a by c. We say $a_0$ because no one thing is everything that is a, i.e. a's method of constructing what we perceive to be a in the context in which that identity exists. The 0 of that method being the function itself, the object of thought that is a's identity. A name does not encompass or cover all the traits that make up a. Dependent on the subject constructing the context, the same subject may be known by the function of another partition of its traits, under another name different from the one whose context it was introduced. What we know are functions whose application express our context. The objects of those contexts do not refer to traits but to the features of our world, the context in which that function is applicable. Thus, $a_0$ is the only member of a that has no value save for what it labels. As the zero of all functions, null is both a part of every articulation constructed and remains apart from them for as a function whose operation expresses what is the case, it is not what those things are in and of themselves. A function cannot be interchangeable with what it puts to use, e.g. $a_1...a_n$ as there would be contradictions, a is only such because it is not-a. What "a" is known as is only such because it is not-a but a property, a label, affixed to a placing in within a system of categorization constructed external to the context in which a is asserted.

Let us say the initial object that emerges is c, making for the context of its assertion described as C. So, "if a and b, then c" is such that c is the case, given (a×b). C describes some subset of $a_n$ and $b_n$, whereby c is less than or equal to $b^a$. If c=(a×b), we know $a_0$ and $b_0$ remain at the tails of the conjunction, the intersection of pair-wise groupings of a,b. Let's make the least assumptions so that $(a_n \times b_n)=c$, leaving $f(a_0)$ or $f(b_0)$. We have conditionalized on some partition of the domain b to a constructing the context in which c can be asserted. This is had by some operation over a and b.

The domain is set as a≤b and an operation over that domain are pairings of a,b up to b, such that,

(a×b)≤c

The tautology c≤c, when c="c" states that "c"=$f(c_0)$=C($a_i,b_i$)

Now suppose there is an operation over c and a. No problem, for in the domain $b^a$, what is not-a is $f(a_0)$="a". Thus c and "a" produce a c-successor in $b^a$. Chomsky's (2019) issue does not arise for if (a×b)≤c leaves what is not-a or not-b or c, because c given (a×b), this operation over c encounters not-a but a's identity. An identity is a function of the label tagged to the set a, what is known as "a" in c. The function "a" plays in c is cited by that operation. $a_0$ is not a physical trait of a, but "a" is an object of thought. Functions are these objects of thought so that regardless of how an "a" may appear, one knows the use of that term within given contexts such that if that object is present, that term can be applied to some end. The method of constructing c, then, is $f(f(a_0),n)=c_n$. Each c-successor, as output of this function, expresses the concept of "c" in the context appropriate for its assertion. Now well-defined, where and when $f(c)$ is applicable it participates in the construction of that context. That label obtains because the application of that function is valid based on its antecedent applications onward from initial use. This chain of succession is proven valid because the operation articulating c, (a×b), was so and thus provides the subject of that context as the conditions utilized to construct it. Each application, although unprovable in the context of its assertion, can be added as an axiom to a successor such that C=0=$f(c)$ as name is constant and C'=$C_1$, $C_1$'=$C_2$, etc. Each context is considered more complete than the previous assertion according to Turing Ordinal logics.

The identity "a" is not a in total, so this operation cites the function "a" plays in the lexicon of c. Now part of c's endowment, this operation represents an internal merge of the functions within c's lexicon. Reference is to c's endowment whose content is the functions that constructed the conditions in which it emerged. It is from the collection of these functions that c can utilize conditions themselves as a tool to articulate an identity appropriate to context of assertion, "a" or "b". What is not-a or not-b, derived from the antecedent of (a×b)≤c in the domain $c^{(a \times b)}$ such that c, given (a×b), is $a_0$ or $b_0$. These names of a and b are accessible to c which assigns them a function, $f(f(…))$. So with (c×a) there is no violation for (c×a)≤c, therefore c' in the domain of c.

On its own, a name labels no one thing, a function without argument. If applicable everywhere, names determine nothing in particular. Exchanging name for what is named loses us the relation between label and entity that we wished to maintain. Null being a member of a as an object of thought and apart from a for it is the function over null, the set of no one thing, what is null within context, then, is the repertoire of the entities from which what was considered "a" were constructed until taken up by c. Null is the domain populated by all articulated and their function of construction, thus what makes articulations possible as the zero of those functions: $f(0)$ for a,b,c. . .

The objects available are the lexical functions with which we construct the contextual descriptions of our world. As assertion implies context of use, it follows that other subjects are doing the same to and with us, relating our contexts and forming the state of affairs we inhabit because we are mutually implicated in its construction. This is how we acquire new names, appropriate to context but not caused by context, and, therefore, not randomly applied. This explains how functions are passed down for they are shared in the form of life within states of affairs.

E.g. One asks someone for a "bottle"; the other does not hand them the actual bottle first named as such, nor do they expect to be handed the same exact bottle every time; one trusts other was given the functional content of that term so when asked, and that label applies, they "know," regardless of appearance, the thing whose function is its capacity to carry liquid. In other contexts, which would make that term inappropriate, the object is called a vase, the functional content of that term is for flowers, although it remains the same physical thing. One refers to the context in which that capacity's label applies, the function of that term defined, and whose argument is the object of thought. This does not apply to physical things which can always be put to different uses, obtain contexts in which their capacities, uses, are different, thus obtaining different names.

With each subject possessing an operation for the construction of contexts from the finite means available, this operation produces infinite arrangements of functions, resulting in an order of "operations"/applications that produces different realities as constructs from the finite material of actual world. I feel this is how we begin to conceive of dealing with the interface problem, how the activity of our thoughts are actualized in the extra mental world. Mainly, inspired by the nondeterminate but finite things composing the actual world, through this recursive capacity to make infinite use of finite means, we develop a finite lexicon of functions qua objective thoughts. It is implausible that we inhabit every context possible from which to experience the actual world. Through this activity, our thoughts interact with the things of the world as we find them.

Most models constructed by our activities can be proven inconsistent. Not a cause for worry, but motivation to create alternatives. The ones that are valid, most often say very little at all, opening us up to what we do not know. Logic seems to lead to openings, rather than closures. As a famous philosopher once said, "the world is everything that is the case." What is the case are beliefs we have decided are true; facts as named objects of knowledge. However, facts can be arranged to tell a multitude of stories, not all of which are true and mostly undecidable. If it cannot be proven consistent, we can take that unprovable statement as an axiom to construct a successor-context. The assertion of that axiom being the name of that successor-context, seemingly more complete than the last, changing our states of affairs for the better. Just because one identity can be asserted and another can as well it does not necessarily follow that both are valid together. Each context of applicable assertion regarding them separately does not mean a context can be constructed in which they both are applicable. Just the same, if two identities are asserted together, it does not necessarily follow that they can be asserted on their own. Knowledge on its own, the zeroes of functions, do not count for very much until applied; either in contexts in which we understand what is the case, because that context is shared with others who are also constructing it; or applied within a system of statements constructed within that body of knowledge for which a contradiction proves the system constructed invalid. If invalidated, this does not negate what the knowledge therein. Athena's helm made her invisible, yet her wisdom was ready for action, her tools, sharpened, could not be taken away.

A name, whose content alone is its function, if held absolute, as universally applicable, determines no one particular thing. Black-ness is applicable across contexts and thus remains functional, no one thing in any one context and, thus, can be expressed utilizing anything. What is made homogeneous attempts to make Black-identity its object. Thus the use of that term implies some context of assertion in which it applies to the objects with which it defines itself. E.g. "w"=$f(w_0)$ if constant; if held universally homogeneous, w→0 asserts that it is something other than a function, for tis dominance is defined by b's assertion as, 0→b, if no one thing then b, an operation producing the value which constitutes, thereby defines, w's position. Thus, if w="b", then w=$f(b_0)$. Defining White-identity as a universal constant reveals its inconsistency, for the class of all Blacks, B-identity, is a function without argument. Thus "White," as sole object, is only applicable in one context. As such, w wipes itself out for all that is produced is b-ness, i.e. what is not-w. This is not to say that blackness produces "white"-ness, but that White-identity has no means of production on its own; we only see w's violent attempt to lay claim to some entity perpetually on the move. Thus, 1=0; or taking Fanon's notion, Black given White (B/W), since White names Black-identity, 0/1 which is 0. Only blackness survives the operation; White-identity is unsustainable. Homogeneity over no one thing, 1/0, generates an infinite number of issues for the system of categorization. Either White supremacy is inconsistent, or blackness is infinitely generative; the disjoint union of both encompasses our current state. Made interchangeable with those things, it cannot stand on its own, thus applicable in one and only one context if at all and, therefore, is empty. Its assertion cannot justify itself for lack of that justification having a means of justification—dominant over everything, therefore, no one particular thing, which does not say anything about dominance. Blackness, having to produce value across states of affairs, by definition an operation, not a thing, has a means of expression, a subjectivity constructing the constituents of the various categories under which it operates by the properties apparent in the contexts in which it is asserted. From no one thing, in light of historical antecedents, blackness is not empty for full of functional capacities, the objects of thought expressing its subject in those contexts and requires no external justification because it is exhibited by the means of construction.

Names have no value in themselves, but as such create value because they are is exportable, lending themselves to racialized capitalist regimes. Who controls the name, controls the value. However, possessing that name is not valuable itself. What is attributed to that name is valued, its function and, therefore, the contexts to which it's exported as an object for consumption. For those whose form of life *is* this function—animating these terms, indicating the operation whose function expresses the concept of what it is to be in that context—this does not prove harmful. It shows that that function can be put to uses against regimes of control or, in the least, circumvent these systems. Functions themselves are illegible in that system, as zeroes they are present but not represented in context. Only their products appear. As such, representation is not everything. Subjectivity is in excess of categorization. For those which forms of life are rendered a thing, they consume the means of theirs or others subordination.

[iii] APPENDIX III

From this study we can draw a tree diagram of the structure of states of affairs in order to show how this concept comes to fruition. Let's take a (C)ontext created by (ID)entity assertion. With the initial assertion of an ID, the ID is identical to context, for its use constructs the context in which it is applicable. So,

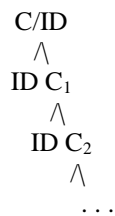

$$\begin{array}{c} C/ID \\ \wedge \\ ID \;\; C_1 \\ \wedge \\ ID \;\; C_2 \\ \wedge \\ \ldots \end{array}$$

C can take another C for the same ID or the initial assertion of another ID=C,

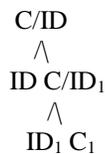

$$\begin{array}{c} C/ID \\ \wedge \\ ID \;\; C/ID_1 \\ \wedge \\ ID_1 \;\; C_1 \end{array}$$

We will specify the concept of this operation as the following. Take a (w)orld as context, a (n)ame as identity, and the (r)elation that name has within context to be that name's function. The function of a name indicates its meaning in and only in the context of that identities assertion. The construction rules above are based on the idea that subjectivity is a vector. Its initial state is no one thing, thus basis for everything. As such, a subject expressed by some operation s, at $s(0)$ is $0=w-n+r$, where $w=n-r(n)$—see Sarrus' rule and Laplace's co-factor expansion for matrix notation $<+, -, +>$. This is the case as $r(0)$ entails that the assertion of a name/ID indicates the world/context of its assertion. From this notation, using the cross product multiplication of vectors we can model how states of affairs emerge by a recursively constructed relation between the contexts of assertion for two or more identities, whose function indicate subjects which may or may not be of the same endowment but, nevertheless, participate in the construction of a state of affairs. This cross product entails the set, here state, in which the disjoint union of the pair-wise coupling of its constituents, i.e. set-recursive multiplication of these well-ordered triples, is the structure composed by the affairs of various subjects. The inner relation of the three aspects of a subject paired with others displays the operation of subjectivity projecting an identity in the context in which it is applicable, because actively constructing that context. So we have,

w

$$\bigwedge_{n\ r(n)}$$

In this way, a name/ID is such that its use indicates the context of its assertion and that identities relation, i.e. function, within and to that context. Therefore the content of the assertion of an identity is its function, not some thing. Thus a name is $<n,r(n)>=<n,<r,n>>$. E.g., an 'I' utterance creates initial context of use and when used again refers to that initial context, thus relating name to that context. The name, a function of that context. As a function of context application is appropriate in context of that line because successors in that line, by that ID utterance, captures the contexts of use as object(s) of that proposition. Subjects considered vectors, vectors having direction and magnitude, display the operation as subjectivity, not some other thing. The context they inhabit is initially determined at its zero,

$0=w-n+r(n)$

So for our tree above, $r(n)=n-w$ standing for an initial state $r(0)$. ID has no scope on its own until applied. This can be simplified so as to state that w, the initial context, is

$w=n-r(n)$

where $w=0$ means that the name or ID of that context is yet to be applied.

When applied, $n=w+r(n)$, which is to say the context indicated by its assertion and it's function constructs that context.

By our diagrams we see how finite contextual material output becomes input for the expression of a successor context in which that ID is appropriate for assertion. Take a context A. The self-reflexive citation of A cannot be proven in A, but if added as an axiom to A creates a set in which A is a member, a successor-context A'. This statement can be proven in that successor context. The context that captures all those contexts represents a separate context and so on. Each addition of a self-reflexive axiom creates a new ordinal in accordance with Turing logics with the hope that each operation makes that context more complete. For example,

$A_n'=A_{n+1}, A_k=B_0$

for all $n<k$, so that $A_0$, $B_0$, etc. are unified in null. If $n=k$, the context marked as such means that the ID, i.e. axiom, is applicable everywhere, therefore nowhere, so that k is not equal to any one n. An initial objection would be that a context can never be determined complete, only work towards completion. This is no problem for us, for to assume that a context is complete makes something which is inherently an operation, subjectivity, interchangeable with that context, a thing, therefore rendering what by definition is not a thing homogeneous and static. If a subject comes to the end of its operations, this marks the end of its life, so what does it care if the context is complete. However, by virtue of our discussion of endowments, if a single subject is made interchangeable with its identity, marking the end of its life, the form of life constructed does not end for subjects utilizing that endowment actively and continue to construct the form of life expressed by their operation over and within that endowment for $A_k=B_0$, thus B' contains A, even though $a_1 \neq b_1$. Considering blackness in this way, one of its identities may be an historical object, but blackness qua subjectivity is an immediate thing, in excess of any one historical frame. Its expressions are appropriate to historical as well as contemporary contexts but not necessarily caused by contextualized reference frames, for its operations actively creates contexts.

Subjectivity, then, is an inner creative capacity for everyone, therefore no one in particular. It creates particular utterances indicative of the individual it represents, dependent upon context and based on the infinite use of finite means appropriate to the context inhabited. The same capacity, different context, entails different output, creating unique individuals based on a restriction due to the finite means to articulate that individual within that context. Nominal-identity is always context dependent. Same subject appears different because contexts are different.

Regarding blackness, the use of any ID indicates the context in which it is appropriate to assert that ID because its use constructs that context. The indication of the operation constructing that context is the subject. Therefore blackness, the function of that operation, cannot be void for there are a multitude of context in which it inhabits and creates under different IDs. The same name can come under different functions. Just because one has access to a context because of their participating in the construction of the state of affairs of which it is a part, does not mean that individuals have access to each other's endowments. Functions being the content of their respective assertions, a function outside of its endowment fails. Different individuals may participate in the construction of same endowment, but only by way of an affinity of operations making for composite functions over IDs created by pair wise set of input/output claims. In this we have the creation of a collective form of life.

[iv] APPENDIX IV

In a 2019 lecture Saul Kripke gave at the University of London, he would state that subject's use their conditions. Those conditions describe the speaker. I say, the operation of that speaker's utilization of a name whose functional content is an assertion of itself is the context of these statements. These can be linked to the previous context of assertion all the way back to the initial baptism, the dubbing of that name. The function of that initial tag marks the occasion of its use; that is when that function obtained the basis for its arguments, which no longer need be present in order to utilize that name. Use refers to the context which obtained that object of thought. Thus, fulfilling the conditions of "referring" making it so recourse is only to the reference frame constructed, not some extra-mental thing that may or may not be there in the long succession of that names use.

Kripke's example of being alone in a room and using the name Gödel provides a great case-study. Description of the traits of Gödel are context dependent. The use of that name does not necessitate that Gödel actually be in the room. The statement describes no extra-mental thing, only the context of its assertion, "the one who formalized the incompleteness theorem." In this empty room there is no way to "know" thereby apply that name to the physical thing that is Gödel. The context of assertion has no thing there save a function relating assertion to context. Within that room, the assertion is appropriate because the functional content of that name is applied in the context of mathematical logic, not a thing. We cannot say that the speaker "knows" the actual person Gödel, nor does the speaker describe any predicament in which they and Gödel were together by the use of that term. The former implies that physical acquaintance is not necessary to use of term, only the possession of a function of the term is needed which also does not entail that the arguments of context and application both obtain that function indicating appropriate use. The latter makes use ambiguous because it refers not to one object but two, Gödel and speaker, revealing that there must be some connection between term and context, not term and thing.

It just so happens to be correct but continued use of this name does not entail access to the citational chain which licenses its use. There is no way to ascertain the context in which its assertion was applicable for the speaker does not share the functional content of that chain because not introduced nor actively participating in the ordinally constructed contexts referring back to the occasion when the one tagged Gödel did what they did at that time. This brings up the point that the functional content is evoked with each utterance of that name up through the citational chain. Thus, the context of assertion is empty, the speaker does not truly "refer" to any one thing save to themselves. Their assertion extends to that and only that context making the identity asserted contingently true, flying in the face of identity theory but proving that the act of referring and identity are not the same thing. So the use of the name obtained no object save its context, however, we can show how, when, and why it would be valid. It is only valid when one is introduced and actively participates in the state of affairs in which they have access to the function relating assertion with its context of application. In looking at names as only having functional content, and how assertion implies context, i.e. the reference frame qua categorical system employed that indexes assertion to the context constructed, we see how systems of reference and their ontological implication, e.g. racial identity, evolve and, when their deployment is proven inconsistent, can change.